\documentclass[12pt]{amsart}

\usepackage{amssymb}
\usepackage{graphicx,color}
\usepackage{amsmath}
\usepackage{cleveref}
\usepackage{paralist}
\usepackage{tikz-cd}
\usepackage{tikz}
\usetikzlibrary{matrix}

\newtheorem{theorem}{Theorem}[section]
\newtheorem{proposition}[theorem]{Proposition}
\newtheorem{lemma}[theorem]{Lemma}
\newtheorem{corollary}[theorem]{Corollary}

\theoremstyle{definition}
\newtheorem{definition}[theorem]{Definition}

\newtheorem{example}[theorem]{Example}
\newtheorem{problem}[theorem]{Problem}
\newtheorem{question}[theorem]{Question}

\newtheorem{remark}[theorem]{Remark}

\newcommand{\ZZ}{ \ensuremath{\mathbb{Z}}}

\newcommand{\C}{\mathcal{C}}

\makeatletter
\def\moverlay{\mathpalette\mov@rlay}
\def\mov@rlay#1#2{\leavevmode\vtop{%
   \baselineskip\z@skip \lineskiplimit-\maxdimen
   \ialign{\hfil$\m@th#1##$\hfil\cr#2\crcr}}}
\newcommand{\charfusion}[3][\mathord]{
    #1{\ifx#1\mathop\vphantom{#2}\fi
        \mathpalette\mov@rlay{#2\cr#3}
      }
    \ifx#1\mathop\expandafter\displaylimits\fi}
\makeatother

\newcommand{\lk}{{\mathrm{lk}}}
\newcommand{\st}{\mathrm{st}}
\newcommand{\BC}{\mathrm{BC}}

\newcommand\bst{\mathrel{\overset{\makebox[0pt]{\mbox{\normalfont\tiny\sffamily bst}}}{\approx}}}
\newcommand\crs{\mathrel{\overset{\makebox[0pt]{\mbox{\normalfont\tiny\sffamily crs}}}{\approx}}}
\newcommand\sh{\mathrel{\overset{\makebox[0pt]{\mbox{\tiny\sffamily sh}}}{\approx}}}
\newcommand\bsh{\mathrel{\overset{\makebox[0pt]{\mbox{\normalfont\tiny\sffamily bsh}}}{\approx}}}
\newcommand\sharr{\mathrel{\overset{\makebox[0pt]{\mbox{\normalfont\tiny\sffamily sh}}}{\longmapsto}}}

\numberwithin{equation}{section}

\begin{document}

\title{Balanced shellings and moves on balanced manifolds}

\author[M. Juhnke-Kubitzke]{Martina Juhnke-Kubitzke}
\email{juhnke-kubitzke@uni-osnabrueck.de}
\author[L. Venturello]{Lorenzo Venturello}
\email{lorenzo.venturello@uni-osnabrueck.de}
\address{
Universit\"{a}t Osnabr\"{u}ck,
Fakult\"{a}t f\"{u}r Mathematik,
Albrechtstra\ss e 28a,
49076 Osnabr\"{u}ck, GERMANY
}

\date{\today}

\thanks{
Both authors were supported by the German Research Council DFG GRK-1916.
}
\keywords{simplicial complex, balancedness, shellability, combinatorial manifold, cross-flips}
\subjclass[2010]{05E45, 57Q15}

\begin{abstract}
A classical result by Pachner states that two $d$-dimensional combinatorial manifolds with boundary are PL homeomorphic if and only they can be connected by a sequence of shellings and inverse shellings. We prove that for balanced, i.e., properly $(d+1)$-colored, manifolds such a sequence can be chosen such that balancedness is preserved in each step. As a key ingredient we establish that any two balanced PL homeomorphic combinatorial manifolds with the same boundary are connected by a sequence of basic cross-flips, as was shown recently by Izmestiev, Klee and Novik for balanced manifolds without boundary. Moreover, we enumerate combinatorially different basic cross-flips and show that roughly half of these suffice to relate any two PL homeomorphic manifolds.

\end{abstract}

\maketitle


\section{Introduction}
An \emph{(elementary) shelling} is the removal of a facet $F$ from a simplicial complex $\Delta$ with the additional requirement that the set
$$
\{G\subseteq F~:~G\notin \Delta\setminus F\}
$$
has a unique minimal element. A pure $d$-dimensional simplicial complex is \emph{shellable} if it can be reduced to the $d$-simplex by a sequence of shellings. Shellability naturally extends to more general objects, such as polyhedral complexes and simplicial posets, and it has become an important concept not only in topological combinatorics \cite{Bjoerner} and polyhedral theory \cite{Bruggesser:Mani} but also in piecewise linear topology \cite{Bing,Rourke:Sanderson}, algebraic combinatorics and combinatorial commutative algebra \cite{Stanley-greenBook} as well as poset theory \cite{Bjoerner:Wachs1,Bjoerner:Wachs2,Wachs}. Prominent examples of shellable simplicial complexes comprise e.g., triangulations of $2$-spheres \cite{Danaraj:Klee} and boundary complexes of polytopes, as shown by Brugesser and Mani \cite{Bruggesser:Mani}. The latter was used by McMullen in his proof of the Upper Bound Theorem, providing tight upper bounds on the face numbers of convex polytopes \cite{McMullen}. Shellability also places strong conditions on the topology of a simplicial complex since it is well-known that any shellable simplicial complex is homotopy equivalent to a wedge of spheres \cite{Danaraj:Klee74}. In particular, any shellable pseudomanifold is homeomorphic to a ball or a sphere. Unfortunately, the converse is not true, i.e., there exist combinatorial balls and spheres that are non-shellable  (see e.g., \cite{Rudin,Ziegler} and \cite{Hachimori:Ziegler,Lickorish91}). 

Allowing not only shellings but also the inverse operations, Pachner \cite[Theorem 6.3]{Pac1} could show the following result:
\begin{theorem}\label{thm:Pachner}
Two combinatorial manifolds with boundary are PL homeomorphic if and only if they are related by a sequence of shellings and inverse shellings.
\end{theorem}

Our focus in this article lies on balanced simplicial complexes, i.e., pure simplicial complexes whose underlying graph admits a ``minimal'' coloring. Those complexes were originally introduced by Stanley \cite{St79} and examples include barycentric subdivisions of regular CW complexes, Coxeter complexes and Tits buildings. While, clearly, shellings preserve balancedness, inverse shellings might add a new possibly monochromatic  edge, which destroys balancedness. In particular, if $\Delta$ and $\Gamma$ are balanced PL homeomorphic manifolds, the sequence provided by \Cref{thm:Pachner} might contain non-balanced simplicial complexes in intermediate steps. Our main result shows that this obstruction can be avoided.

\begin{theorem}\label{mainresult}
Two balanced combinatorial manifolds $\Delta$ and $\Gamma$ with boundary are PL homeomorphic if and only if they can be connected by a sequence of shellings and inverse shellings that preserves balancedness in each step.
\end{theorem}
\Cref{mainresult} provides a positive answer to Problem 1 in \cite{IKN}, posed by Izmestiev, Klee and Novik. Our proof technique combines ideas of Pachner's proof of \Cref{thm:Pachner} and methods developed and employed in \cite{IKN}. As a key step, we use those ideas together with a result by Casali \cite[Proposition 4]{Cas} to show the following:

\begin{theorem}\label{analogIKNBoundary}
Let $\Delta$ and $\Gamma$ be balanced combinatorial manifolds with $\partial\Delta\cong\partial\Gamma$. Assume moreover that the isomorphism preserves the coloring. Then $\Delta$ and $\Gamma$ are PL homeomorphic if and only if they are related by a sequence of cross-flips.
\end{theorem}
The previous result provides an analog of Theorem 1.2 in \cite{IKN}, where the corresponding statement was shown for closed manifolds. Roughly speaking a cross-flip is the balanced analog of a bistellar flip and it substitutes a subcomplex of the boundary of the cross-polytope with its complement. (We defer more details and the precise definitions to \Cref{sect:background}.) With \Cref{analogIKNBoundary} in hand, the strategy to show \Cref{mainresult} is to first reduce to the situation that $\Delta$ and $\Gamma$ have the same boundary and then to convert each cross-flip needed to transform $\Delta$ into $\Gamma$ into a sequence of shellings and inverses. The latter requires two ingredients: First, the construction of shellings for particular subcomplexes of the boundary of the cross-polytope, relative to their boundaries (\Cref{theorem:shelling}), and second, building a ``collar'' around a manifold in order to protect its boundary~--~an idea which goes back to Pachner (\Cref{collar}).
 
It was shown in \cite{IKN} that in order to relate any two closed balanced PL homeomorphic manifolds it is enough to consider a restricted set of moves, referred to  as \emph{basic} cross-flips (see \Cref{Sect:cross-flips} for the precise definition). \cite[Problem 2]{IKN} asks for a description and the number of combinatorially distinct basic cross-flips. We provide an answer to this question, which can be summarized as follows (see Theorems \ref{theorem:NumberCrossflips} and \ref{thm:reducible} for the detailed statements):

\begin{theorem}\label{thm:BasicCrossflipsSummarized}
There are $2^{d+1}-1$ combinatorially distinct basic cross-flips in dimension $d$, out of which $2^d$ are sufficient to relate any two $d$-dimensional PL homeomorphic balanced manifolds without boundary or with the same boundary.
\end{theorem}

The enumeration of  combinatorially distinct cross-flips relies on a detailed study of their combinatorics. For the proof of the second part of \Cref{thm:BasicCrossflipsSummarized} we construct a set $M$ of basic cross-flips with $|M|=2^d$ such that any other basic cross-flip can be expressed as a combination of cross-flips in $M$.

The layout of this article is as follows. \Cref{sect:background} provides necessary background on simplicial complexes and the combinatorics of local moves. \Cref{Section:balancedShelling} contains the proof of \Cref{mainresult}. In \Cref{sect:BasicCrossflips} we prove the statements of \Cref{thm:BasicCrossflipsSummarized}. We end this article with some open problems in \Cref{sect:openProblems}.

\section{Combinatorics of and moves on simplicial complexes}\label{sect:background}
In this section, we provide background on the combinatorics of simplicial complexes and simplicial moves (including stellar subdivisions, bistellar flips and cross-flips). The last part of this section discusses shellability and (inverse) shellings. 

\subsection{Simplicial complexes and combinatorial manifolds}
We start with several definitions. 
For a detailed exposition of this subject we refer to Stanley's book \cite{Stanley-greenBook}. 
An (abstract) \emph{simplicial complex} $\Delta$ on a (finite) vertex set $V=V(\Delta)$ is a collection of subsets of $V(\Delta)$ that is closed under inclusion. Elements $F\in\Delta$ are called \emph{faces} of $\Delta$, and maximal faces (with respect to inclusion) are called \emph{facets}. We use $\mathcal{F}(\Delta)$ to denote the set of facets of $\Delta$. The \emph{dimension} of a face $F$ is $\dim(F)=\left|F \right|-1$ and the \emph{dimension} of $\Delta$ is $\dim(\Delta)=\max\left\lbrace \dim(F)~:~F\in\Delta\right\rbrace$. $0$-dimensional and $1$-dimensional faces are also called \emph{vertices} and \emph{edges} of $\Delta$, respectively.  If all facets of $\Delta$ are of the same dimension, then $\Delta$ is called \emph{pure}. 
 A \emph{subcomplex} $\Gamma$ of $\Delta$ is any simplicial complex $\Gamma\subseteq \Delta$. We call such a subcomplex $\Gamma\subseteq\Delta$ \emph{induced} if any $F\subseteq V(\Gamma)$ with $F\in\Delta$ is a face of $\Gamma$. Given a collection of faces $F_1,\dots,F_m$ of a simplicial complex $\Delta$, we denote by $\left\langle F_1,\dots,F_m\right\rangle$ the smallest simplicial complex that contains $F_1,\ldots,F_m$, i.e.,
$$
\left\langle F_1,\dots,F_m\right\rangle=\{F\subseteq V(\Delta)~:~F\subseteq F_i\mbox{ for some } 1\leq i\leq m\}.
$$
 Given a face $F\in\Delta$, the \emph{link} and the \emph{star} of $F$ in $\Delta$ are two simplicial complexes providing a local description of $\Delta$ around $F$: 
\begin{equation*}
\lk_{\Delta}(F)=\left\lbrace G\in\Delta~:~ F\cup G\in\Delta, F\cap G=\emptyset\right\rbrace \; \mbox{and} \;
\st_{\Delta}(F)=\left\lbrace G\in\Delta~:~ F\cup G\in\Delta\right\rbrace.
\end{equation*}
The \emph{deletion} of a face $F\in \Delta$ describes the simplicial complex $\Delta$ outside of $F$:
$$
\Delta\setminus F=\left\lbrace G\in\Delta~:~ F\nsubseteq G\right\rbrace. 
$$
Similarly, the \emph{deletion} $\Delta\setminus \Gamma$ of a subcomplex $\Gamma$ from a simplicial complex $\Delta$ is defined as
$$
\Delta\setminus\Gamma=\left\langle \mathcal{F}(\Delta)\setminus\mathcal{F}(\Gamma)\right\rangle.
$$ 
The \emph{join} of the simplicial complexes $\Delta$ and $\Gamma$ on disjoint vertex sets is
$$\Delta*\Gamma=\left\lbrace F\cup G~:~F\in\Delta,G\in\Gamma\right\rbrace.$$ 

A \emph{combinatorial $d$-sphere} (or \emph{PL sphere}) is a simplicial complex that is PL homeomorphic to the boundary of a $(d+1)$-simplex $\partial\sigma^{d+1}$. Similarly, a \emph{combinatorial $d$-ball} is a simplicial complex that is PL homeomorphic to a $d$-simplex $\sigma^d$.	Here and in the following, we use $\sigma^d$ and $\partial\sigma^d$ to denote a $d$-simplex and its boundary complex, when the vertex set is not important. 

A $d$-dimensional simplicial complex $\Delta$ is a combinatorial $d$-sphere if and only if $\Delta$ and $\partial \sigma^{d+1}$ have a common subdivision. Moreover, any combinatorial $d$-sphere is homeomorphic to a $d$-sphere and for $d\leq 3$ also the converse is true. For $d\geq 5$ however, there exist $d$-dimensional simplicial complexes that are homemorphic to a $d$-sphere but that are not combinatorial spheres. An instance for this phenomenon is provided by the $(d-3)$-suspension of the Poincar\'{e} $3$-sphere ($d\geq 5$), which is homeomorphic to a $d$-sphere \cite{Cannon,Edwards}, but which is not a combinatorial sphere, since the Poincar\'{e} $3$-sphere occurs as one if the links (see also \cite{Bjoerner:Lutz} for examples of topological spheres with few vertices that are non-PL spheres). The problem of whether any $4$-dimensional simplicial complex that is homemorphic to a sphere is also a combinatorial $4$-sphere is open.

	A closed \emph{combinatorial $d$-manifold} is a connected simplicial complex $\Delta$ such that for every non-empty face $F\in\Delta$ the link $\lk_\Delta(F)$ is a combinatorial $(d-|F|)$-sphere. A \emph{combinatorial $d$-manifold with boundary} is a connected  simplicial complex $\Delta$ such that for every non-empty face $F\in \Delta$ the link $\lk_\Delta (F)$ is either a combinatorial $(d-|F|)$-sphere or a combinatorial $(d-|F|)$-ball. 
If $\Delta$ is a combinatorial $d$-manifold with boundary, then its \emph{boundary complex} $\partial\Delta$ is defined as
$$
\partial\Delta=\left\lbrace F\in\Delta~:~ \lk_{\Delta}(F) \text{ is a combinatorial ball} \right\rbrace\cup\{\emptyset\}.
$$ 
Equivalently, the boundary complex of $\Delta$, which is itself a closed combinatorial $(d-1)$-manifold, is the simplicial complex whose facets are the $(d-1)$-faces of $\Delta$ contained in exactly one facet of $\Delta$. We will take this as the definition of the boundary complex of any simplicial complex $\Delta$. So, in particular $\partial\Delta=\emptyset $, if $\Delta$ is a closed manifold. If $F\in \Delta$ is a face of a simplicial complex $\Delta$, we write $\partial F$ for its boundary complex, i.e., $\partial F=\{G\subsetneq F\}$. 
We define the \emph{interior}  $\overset{\circ}{\Delta}$ of a combinatorial $d$-manifold $\Delta$ as the set of  all faces of $\Delta$ not in the boundary complex. 

An easy but important combinatorial invariant associated to a $d$-dimensional simplicial complex $\Delta$ is its \emph{$f$-vector} $f(\Delta)=(f_{-1}(\Delta)\dots,f_d(\Delta))$, where $f_{j}(\Delta)$ denotes the number of $j$-dimensional faces of $\Delta$. Often, it is more convenient to consider the so-called \emph{$h$-vector} $h(\Delta)=(h_0(\Delta),\dots,h_{d+1}(\Delta))$ of $\Delta$, defined by 
$$h_j(\Delta)=\sum_{i=0}^{j}\binom{d+1-i}{d+1-j}f_{i-1}(\Delta) \quad \mbox{ for }\quad 0\leq j\leq d+1 .$$

\subsection{Stellar moves and bistellar flips}
In this section we define different local moves on simplicial complexes and state well known results on the equivalence classes determined by such moves.

Given a simplicial complex $\Delta$ and a face $F\in \Delta$, the \emph{stellar subdivision} of $\Delta$ at $F$ is the simplicial complex 
	$$\text{sd}_{F}(\Delta)=\Delta\setminus F\cup (\langle v\rangle *\partial F*\lk_{\Delta}(F))$$
	where $v\notin\Delta$ is a new vertex. 
If two simplicial complexes can be transformed one into the other by a sequence of stellar subdivisions and their inverses (stellar welds), we say that they are \emph{stellarly equivalent}. Clearly, neither  subdivisions nor welds do affect the topology of a simplicial complex. Indeed the following classical result was shown by Alexander:
\begin{theorem}
	\cite[Theorem 10:3]{Ale}\label{AlexanderBig}
	Two simplicial complexes $\Delta$ and $\Gamma$  are PL homeomorphic if and only if
	they are stellarly equivalent. 
\end{theorem}	
Several other results in the same flavor exist, e.g., Alexander \cite{Ale} and Newman \cite{Newman} independently showed that edge subdivisions and welds suffice, and Ludwig and Reitzner provided a ``geometric'' version of this result for polytopes \cite{Ludwig:Reitzner}. Moreover, Lutz and Nevo \cite{Lutz:Nevo} proved that PL homeomorphic flag simplicial complexes can be transformed into each other by a sequence of edge subdivisions and welds such that flagness is preserved in each step.  

As the number of facets added by a stellar subdivision at a face $F\in \Delta$ depends on the combinatorics of the link $\lk_{\Delta}(F)$, unfortunately there are infinitely many combinatorially different stellar subdivisions even if both the dimension of $\Delta$ and the dimension of $F$ are fixed. The following set of moves, that was introduced by Pachner \cite{Pac3}, remedies this situation by providing finitely many moves for each dimension.

Let $\Delta$ be a $d$-dimensional simplicial complex and assume that there exists a face $A\in\Delta$ such that $\lk_{\Delta}(A)=\partial B$, for some $B\notin\Delta$. A \emph{bistellar flip} (or \emph{bistellar move}) on $\Delta$ is the operation $\chi_{A,B}$ defined by
	$$\Delta\mapsto \chi_{A,B}(\Delta)=\Delta\setminus(\left\langle A\right\rangle *\partial B)\cup(\partial A*\left\langle B\right\rangle ),$$
		i.e., a bistellar flip exchanges $\left\langle A\right\rangle *\partial B$ with $\partial A*\left\langle B\right\rangle$. Clearly, the inverse of a bistellar flip $\chi_{A,B}$ is given by the bistellar flip $\chi_{B,A}$. Two simplicial complexes $\Delta$ and $\Gamma$ are called \emph{bistellar equivalent} if they are related by a sequence of bistellar flips. We write $\Delta\bst\Gamma$ in this case.
    Bistellar moves admit a nice and simple geometric description: Indeed, the bistellar flip $\chi_{A,B}$ just replaces the subcomplex $\langle A\rangle* \partial B$ that is isomorphic to a subcomplex $D$ of $\partial\sigma^{d+1}$ which is a $d$-ball, with the complex $\partial A*\left\langle B\right\rangle$ that is isomorphic to the complement of $D$ in $\partial\sigma^{d+1}$. 
		
\begin{figure}[h]
	\centering
	\includegraphics[scale=0.7]{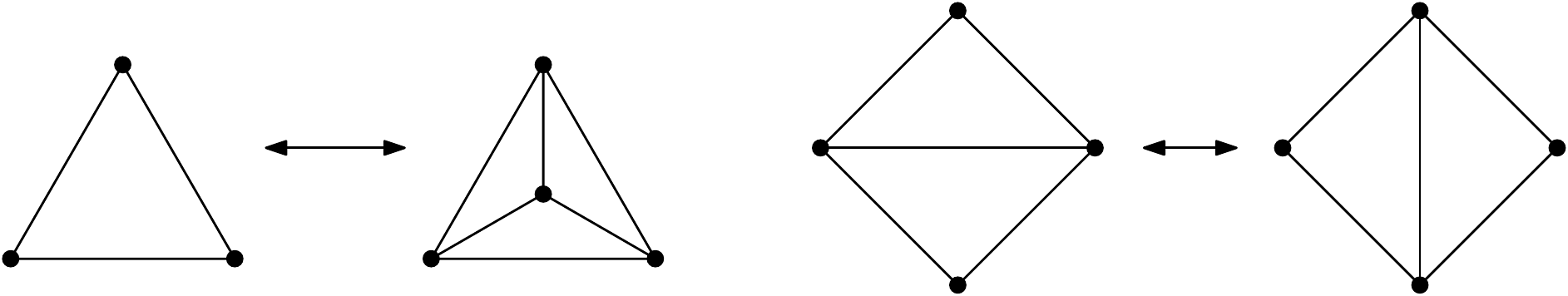}
	\caption{Bistellar flips for $d=2$.}
	\label{bistellar}
\end{figure}

    As $\partial\sigma^{d+1}$ has exactly $d+1$ combinatorially different pure $d$-dimensional subcomplexes (that are all $d$-balls), there are exactly $d+1$ distinct bistellar flips in dimension $d$. Figure \ref{bistellar} depicts all bistellar flips in dimension $2$. The following analog of \Cref{AlexanderBig} is due to Pachner (see also \cite{Lickorish} for a proof).
\begin{theorem}\cite[Theorem 5.5]{Pac1}
	\label{pac_bist}
	Two closed combinatorial manifolds $\Delta$ and $\Gamma$ are PL homeomorphic if and only if $\Delta\bst\Gamma$.
\end{theorem}
Since any bistellar flip can be written as a composition of a stellar subdivision and a weld, the ``If''-part of \Cref{pac_bist} is a direct consequence of \Cref{AlexanderBig}. However, for closed manifolds stellar and bistellar equivalence even turn out to be equally strong. Indeed, Pachner \cite[Lemma 4.8]{Pac1} showed that every stellar subdivision at a face $F$ in the interior of a simplicial complex $\Delta$ such that $\lk_\Delta(F)$ is shellable can be realized by a sequence of bistellar flips. Casali \cite[Proposition 4]{Cas} could improve on this result by proving that the shellability assumption is not necessary. As an almost immediate consequence, she obtained the following analog of \Cref{pac_bist} for manifolds with boundary.

\begin{proposition}\cite[Main Theorem]{Cas}
	\label{same_bound_bist_equi}
	Let $\Delta$ and $\Gamma$ be combinatorial $d$-manifolds with isomorphic boundaries, i.e.,  $\partial\Delta\cong\partial\Gamma$. Then $\Delta$ and $\Gamma$ are PL homeomorphic if and only they are bistellar equivalent. 
\end{proposition}
Observe that, since bistellar flips do not affect the boundary, all manifolds constructed from the sequence of bistellar flips, guaranteed by the previous proposition, have the same boundary. We will use \Cref{same_bound_bist_equi} together with this simple observation in \Cref{Sect:pseudocobordism}. 

\subsection{Bistellar moves and cross-flips on balanced complexes}\label{Sect:cross-flips}
A natural question that arises from the previous sections is, whether there are analogs of \Cref{AlexanderBig} and \Cref{pac_bist} for special classes of simplicial complexes. In the following, our focus lies on balanced complexes whose definition we now recall. 

A simplicial complex $\Delta$ is called \emph{properly $m$-colorable} if there exists a map (a \emph{coloring}) $\kappa:\; V(\Delta)\to \{0,1,\ldots,m-1\}$ such that there are no monochromatic edges, i.e., $\kappa(u)\neq \kappa(v)$ for all $\{u,v\}\in \Delta$. Clearly, a proper coloring of a $d$-dimensional simplicial complex $\Delta$ requires at least $d+1$ different colors. If such a minimal coloring exists, we say that $\Delta$ is \emph{balanced}. Examples include Coxeter complexes, Tits buildings and barycentric subdivisions of regular CW complexes, which in particular means that if a topological manifold admits a (finite) triangulation, then it also admits a balanced one. Note that it follows from the refutation of the Triangulation Conjecture by Manolescu, that for $d\geq 5$ there exist compact topological manifolds of dimension $d$ that are not homeomorphic to a finite simplicial complex \cite{Manolescu}. 
Balanced simplicial complexes have been studied intensively during the last years and many classical results have been proven to exhibit balanced analogs \cite{IKN,Juhnke:Murai,Juhnke:Murai:Novik:Sawaske,KN}.

It is easy to see that stellar subdivisions might destroy balancedness of a simplicial complex. In fact, even if the resulting complex after applying such a move is balanced, in general the vertex colors are not preserved. As the balanced analog of stellar subdivisions, Fisk \cite{Fisk1977,Fisk} introduced balanced stellar subdivisions (see also \cite{Izmestiev:Joswig}). Very recently, Murai and Suzuki \cite{Murai:Suzuki} showed that even in dimension 2 \Cref{AlexanderBig} does not have a balanced analog. 

At first glance, the situation appears similar for bistellar subdivisions. In general, balancedness is not maintained and even if it is, the vertex colors might change. 
 Nevertheless, very recently, Izmestiev, Klee and Novik succeeded in proving the following colored version of \Cref{pac_bist}.
\begin{theorem}
	\cite[Theorem 1.1]{IKN} \label{ikn_bist}
	Let $\Delta$ and $\Gamma$ be closed combinatorial $d$-manifolds that are PL homeomorphic. Assume
	that $\Delta$ and $\Gamma$ are properly $m$-colored, $m\geq d+2$. Then there exists a sequence of bistellar flips that
	transforms $\Delta$ into $\Gamma$ such that each intermediate complex is properly $m$-colored and the flips preserve the vertex colors.	
\end{theorem}

Note that the last result does not cover the case of balanced combinatorial manifolds and indeed, as remarked before, balancedness might be destroyed by bistellar flips. The proof of \Cref{ikn_bist} makes use of a so-called $m$-colorable \emph{pseudo-cobordism} that connects $\Delta$ and $\Gamma$ by a sequence of shellings and inverse shellings. The latter sequence is turned into a sequence of bistellar flips between $\Delta$ and $\Gamma$. However, the coloring of the pseudo-cobordism requires at least $d+2$ colors. A similar idea will be used in \Cref{Sect:ProofMainResult} to prove our main result \Cref{mainresult}.

For balanced complexes the right analog of bistellar flips are so-called \emph{cross-flips}, introduced in \cite[Definition 2.6]{IKN}. Recall that a bistellar flip can be defined by substituting a $d$-ball in $\partial\sigma^{d+1}$ by its complement. For balanced complexes it has turned out that the boundary of the $(d+1)$-dimensional cross-polytope plays the same role as $\partial\sigma^{d+1}$ does for arbitrary simplicial complexes (see e.g., \cite{IKN,Juhnke:Murai, Juhnke:Murai:Novik:Sawaske,KN}). The definition of cross-flips combines those two insights. We make this now more precise. The boundary complex $\C_d$ of the $(d+1)$-dimensional cross-polytope is defined to be the $d$-dimensional simplicial complex on vertex set $V(\C_d)=\{x_i,y_i~:~0\leq i\leq d\}$ and whose facets are all sets $F\subseteq V(\C_d)$ with $|F\cap\{x_i,y_i\}|=1$ for all $0\leq i\leq d$.  
Note that $\C_d$ is balanced, since setting $\kappa(x_i)=\kappa(y_i)=i$, for $0\leq i\leq d$ defines a proper $(d+1)$-coloring and it has turned out to be the minimal (with respect to face numbers) balanced $d$-sphere. We call a pure subcomplex $D\subseteq\mathcal{C}_{d}$  \emph{co-shellable} if $\mathcal{C}_{d}\setminus D$ is a shellable simplicial complex. For the definition of shellability we refer to \Cref{SectionShellability}. 
\begin{definition}\label{Def:cross-flips}
	Let $\Delta$ be a balanced $d$-dimensional simplicial complex and let $D\subsetneq\Delta$ be an induced subcomplex that is isomorphic to a shellable and co-shellable subcomplex of $\mathcal{C}_{d}$. The operation $\chi^*_D$ given by
	$$
	\Delta\mapsto\chi^*_D(\Delta)=(\Delta\setminus D) \cup (\C_{d}\setminus D)$$
	is called a \emph{cross-flip} on $\Delta$. If two balanced simplicial complexes $\Delta$ and $\Gamma$ are connected by a sequence of cross-flips, we write $\Delta\crs\Gamma$.
\end{definition}
The fact that $D$ is shellable and co-shellable directly implies that a cross-flip exchanges a $d$-ball by another $d$-ball sharing the same boundary. In particular, $\Delta$ and $\chi^*_D(\Delta)$ are PL homeomorphic. It is easy to see that $\chi^*_D(\Delta)$ is balanced and that the coloring is preserved. Moreover, the inverse of the cross-flip $\chi^*_D$ is given by the cross-flip $\chi^*_{\C_d\setminus D}$, which justifies the notation $\Delta\crs \Gamma$.  
\begin{figure}[h]
	\centering
	\includegraphics[scale=0.7]{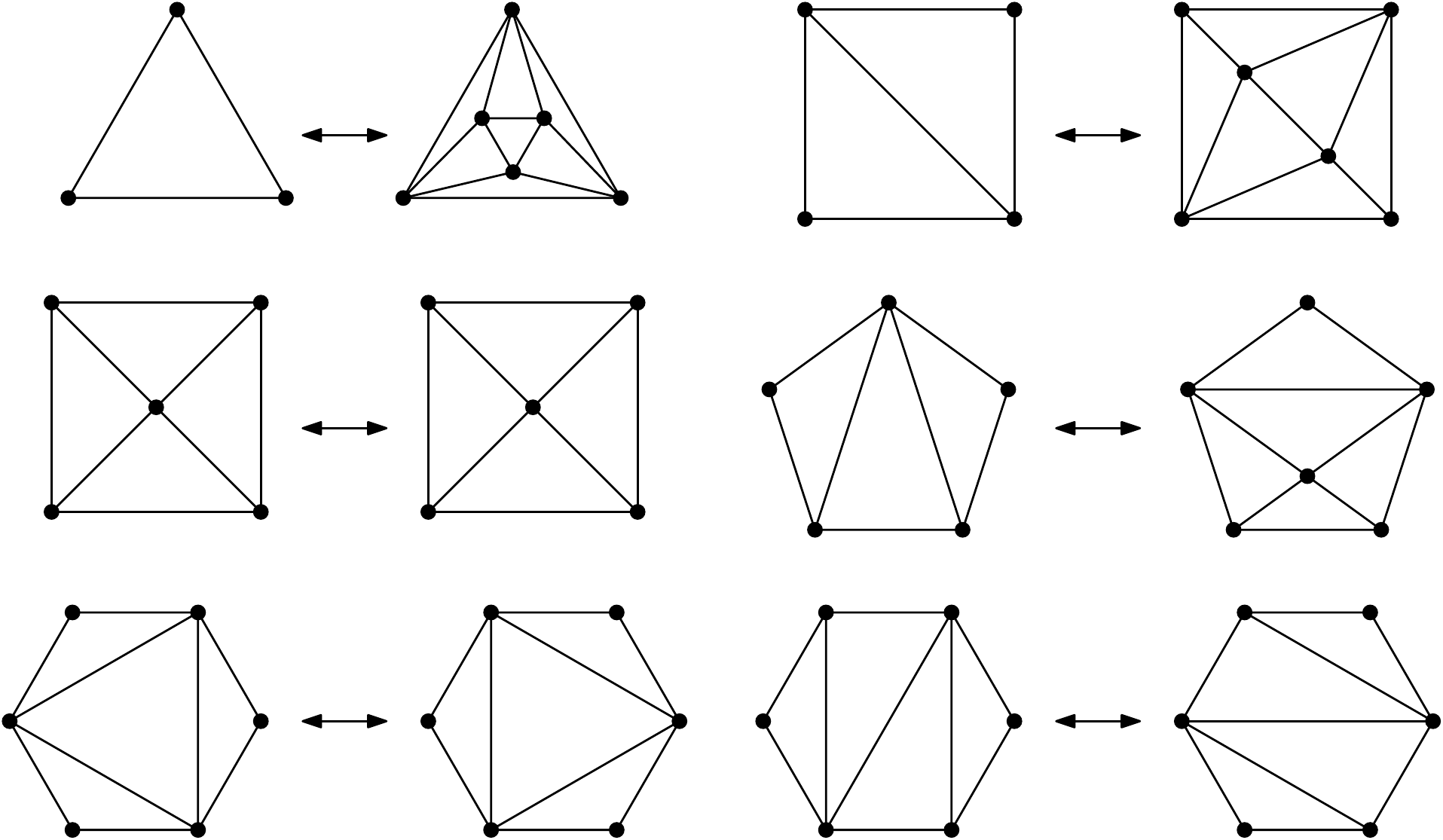}
	\caption{Cross-flips for $d=2$.}
	\label{cross}
\end{figure} 
It is important to underline that the shellability of $D$ implies its co-shellability if and only if $\C_d$ is extendably shellable. For $d \geq 11$, Hall showed $\C_d$ is not extendably shellable \cite{Hall}.

The following result provides the balanced analog of Pachner's \Cref{pac_bist}.
\begin{theorem}\cite[Theorem 1.2]{IKN}
\label{ikn cross}
	Two closed balanced combinatorial manifolds $\Delta$ and $\Gamma$ are PL homeomorphic if and only if $\Delta\crs\Gamma$.	
\end{theorem}
Clearly, there are only finitely many $d$-dimensional cross-flips. However, compared to the number of $d$-dimensional bistellar flips, their number is considerate. It is therefore natural to ask, if \Cref{ikn cross} can be improved by showing that a particular subset of cross-flips suffices. Indeed, in \cite{IKN} it is shown that such a set is provided by the so-called \emph{basic cross-flips}. We now recall this construction, which uses the so-called \emph{diamond operation}:

 For simplicity, we assume that the vertices of the $(d+1)$-simplex $\sigma^{d+1}$ are labelled by $0,\ldots,d+1$. Let $\Gamma\subseteq \partial\sigma^{d+1}$ be a pure and $d$-dimensional subcomplex. Following \cite[Section 3.3]{IKN}, we define another combinatorial $d$-ball $\Diamond(\Gamma)$ that is a subcomplex of $\C_d$ in the following way: For $0\leq i\leq d$: If $F_i=\{i+1,\ldots,d+1\}$ is a face of $\Gamma$, then recursively stellar subdivide $\Gamma$ at $F_i$ and label the newly introduced vertex by $v_{i}$. In particular, the vertex $F_d=\{d+1\}$ gets renamed with $v_d$ in this procedure.  
Note that  $\Diamond(\partial\sigma^{d+1})$ gives the boundary $\C_d$ of the $(d+1)$-dimensional cross-polytope on vertex set $\{0,\ldots,d\}\cup\{v_0,\ldots,v_d\}$ (see \cite[Lemma 3.6]{IKN}). As stellar subdivisions preserve shellability, $\Diamond(\Gamma)$ is a shellable and co-shellable $d$-ball if $\Gamma\subsetneq \partial \sigma^{d+1}$ is a $d$-ball. Those are exactly the balls, one considers in the definition of basic cross-flips.
\begin{definition}\label{def:basic cross-flips}
Let $\Gamma\subsetneq \partial\sigma^{d+1}$ be a $d$-ball. The cross-flip $\chi^*_{\Diamond(\Gamma)}$ is called a \emph{basic cross-flip}.
\end{definition}
The inverse of the basic cross-flip $\chi^*_{\Diamond(\Gamma)}$ is again a basic cross-flip, namely $\chi^*_{\C_d\setminus\Diamond(\Gamma)}$.  
We also want to remark that the diamond operation can easily be extended to any pure balanced $(d+1)$-dimensional simplicial complex $\Delta$ with a specified coloring by first interpreting the vertex colors as vertex labels and then applying the diamond operation to the boundary of every $(d+1)$-simplex. In this way, a balanced simplicial complex can be converted into a cross-polytopal complex. We will use this idea from \cite{IKN} in \Cref{Section:balancedShelling}. 

In \cite{IKN} (essentially the proof of Theorem 3.10; cf., \cite[Remark 3.12]{IKN}), the following improvement of \Cref{ikn cross} was proven:
\begin{theorem}\label{strong thm:crossflips}
Two closed balanced combinatorial manifolds are PL homeomorphic if and only if they can be obtained from each other by a sequence of basic cross-flip.
\end{theorem}

Though the number of basic cross-flip is already considerably smaller than the number of cross-flips, a priori there are still about $2^{d+2}-2$ many (one for each proper subcomplex of $\sigma^{d+1}$ that is a $d$-ball). However, not all of those are combinatorially different. In dimension $2$ there are $11$ combinatorially distinct cross-flips (see \Cref{cross}), but only the $7$ in the first two lines of \Cref{cross} turn out to be basic. Moreover, the move in the middle row on the left is the trivial move, which does not have any effect 

It is hence natural to raise the following problem:
\begin{problem}\cite[Problem 2]{IKN}\label{problem:cross-flips}
 Give an explicit description of basic cross-flips. How many combinatorially distinct basic cross-flips are there? 
\end{problem}
\Cref{theorem:NumberCrossflips} will provide a solution to this question.

\subsection{Shellings and their inverses}\label{SectionShellability}
It is worth remarking that all operations considered so far leave the boundary of a simplicial complex unchanged. Hence, if one wants to connect PL homeomorphic manifolds with boundary another set of moves is needed. This set is provided by shellings and inverse shellings. 

A pure $d$-dimensional simplicial complex $\Delta$ is \emph{shellable} if there exists an ordering $F_1,\dots,F_m$ of the facets of $\Delta$ such that for every $1\leq i\leq m$ the set 
$$\left\lbrace G\subseteq F_i~:~ G\nsubseteq F_j \text{ for }1\leq j\leq i-1 \right\rbrace$$
has a unique minimal element $R(F_i)$, the so-called \emph{restriction face} of $F_i$. The ordering $F_1,\ldots,F_m$ is called a \emph{shelling} of $\Delta$.

Shellings can be used to compute the $h$-vector of a simplicial complex, since one has the relation
\begin{equation}\label{eq:shelling}
h_i(\Delta)=\left| \left\lbrace 1\leq j\leq m: \left| R(F_j)\right|=i  \right\rbrace \right| \mbox{ for } 0\leq i\leq d+1
\end{equation}
(see e.g., \cite[Theorem 8.19]{ZieglerBook}). We will use this to compute the $h$-vectors of diamond complexes in \Cref{sect:numberCrossflips}. 

\begin{definition}\label{definition:shelling}
	 Let $\Delta$ be a pure $d$-dimensional simplicial complex and let $F\in \Delta$ be a facet. Assume that $F$ can be written as $F=A\cup R$, where 
	\begin{enumerate}
	\item[(1)] $\dim A\geq 0$, $\dim R\geq 0$,
	\item[(2)] $A\in \overset{\circ}{\Delta}$.
	\item[(3)] $\partial A *\left\langle R\right\rangle\subseteq\partial\Delta$.  
	\end{enumerate}
	The operation
	$$
	\Delta\sharr\Delta\setminus F
	$$
	is called an (elementary) \emph{shelling} on $\Delta$. The inverse operation is referred to as \emph{inverse shelling}. 
\end{definition}
If two pure simplicial complexes $\Delta$ and $\Gamma$ are related by a sequence of shellings and inverse shellings, we write $\Delta\sh\Gamma$. 

Now let $F_1,\ldots,F_m$ be a shelling of a pure $d$-dimensional simplicial complex $\Delta$ and let us define $\Delta_i=\langle F_1,\ldots,F_i\rangle$ for $1\leq i\leq m$. It directly follows from the definition of a shelling that for every facet $F_i$, we either have $F_i=R(F_i)$, or $F_i$ can be decomposed as in \Cref{definition:shelling} (with $\Delta_i$ in place of $\Delta$). In particular, it follows that shellable combinatorial manifolds with boundary are exactly those simplicial complexes that can be transformed into a simplex by a sequence of shellings (without inverses), which implies that any shellable combinatorial manifold with boundary is a combinatorial ball. Similarly, any shellable combinatorial manifold without boundary is a combinatorial sphere. Once again, shellings and their inverses preserve the PL homeomorphism type and Pachner showed that also the converse is true.
\begin{theorem}\cite[Theorem 6.3]{Pac1}
	\label{pac_shell}
	Two combinatorial manifolds with boundary $\Delta$ and $\Gamma$ are PL homeomorphic if and only if $\Delta\sh\Gamma$.
\end{theorem}
Following the line of discussion of the previous section it is natural to ask, what happens if in \Cref{pac_shell} one assumes $\Delta$ and $\Gamma$ to be balanced. On the one hand, shellings are rather harmless, since no new edges are created and since the resulting complex is a subcomplex of the starting complex. On the other hand, inverse shellings with $\dim(R)=1$ create new edges and those might be monochromatic. In particular, balancedness is destroyed in this case. This motivates the following question by Izmestiev, Klee and Novik \cite[Problem 1]{IKN}:
\begin{question}\label{question:shellable}
Can any two PL homeomorphic balanced combinatorial manifolds with boundary
be related by a sequence of elementary shellings and inverse shellings, such that balancedness (and the coloring) is preserved in each intermediate step? 
\end{question}
If two balanced combinatorial manifolds with boundary $\Delta$ and $\Gamma$ can be connected by such a sequence, we write $\Delta\bsh\Gamma$. Inverse shellings that preserve balancedness will also be referred to as \emph{balanced inverse shellings} in the following. 
It is not hard to see that \Cref{question:shellable} has an affirmative answer if $\Delta$ and $\Gamma$ are balanced shellable balls, since in this case they can be reduced to the simplex only using shellings. However, it is well-known that there are combinatorial balls that are  non-shellable (see e.g., \cite{Rudin} and \cite{Ziegler}).
Our main result \Cref{mainresult} answers \Cref{question:shellable} in the positive in full generality and thereby provides a balanced analog of \Cref{pac_shell}. In particular, together with \Cref{pac_bist}, \Cref{ikn cross} and \Cref{pac_shell} it completes the picture, telling us which moves are necessary to relate any two PL homeomorphic manifolds (with or without) boundary in the balanced as well as in the non-balanced case.

\section{Balanced shellings for combinatorial manifolds with boundary}\label{Section:balancedShelling}
The aim of this section is to prove our main result \Cref{mainresult}. 
The proof will require several intermediate steps and we start with a brief outline of the proof strategy, which should serve as a golden thread in this section.

Let us assume that $\Delta$ and $\Gamma$ are balanced PL homeomorphic combinatorial manifolds with boundary.
\begin{asparaenum}
\item[{\sf Step 1:}] First, via shellings we convert $\Delta$ into a balanced manifold $\Delta'$ such that $\Delta'$ and $\Gamma$ have isomorphic boundaries with the same coloring. (This is \Cref{bound_same}.)\\
\item[{\sf Step 2:}] It follows from \Cref{same_bound_bist_equi} that $\Delta'$ and $\Gamma$ can be connected by a sequence of bistellar flips. Adapting Theorem 4.8, Lemma 5.2 and Corollary 5.3 from \cite{IKN} to our situation, we encode this sequence of bistellar flips by a shellable pseudo-cobordism $(\Omega, \phi,\psi)$ between $\Delta'$ and $\Gamma$. (This is \Cref{pseudo_disj}.)
\medskip
\item[{\sf Step 3:}] Applying the diamond operation to $\Omega$ yields a cross-polytopal complex. As a result, every bistellar flip is converted into a basic cross-flip. (This is \Cref{cross_int}.)
\medskip
\item[{\sf Step 4:}] The last step consists of converting every cross-flip into a sequence of shellings, followed by a sequence of balanced inverse shellings (see \Cref{theorem:shelling}). This step also requires building a balanced collar around a balanced manifold with boundary, an idea already appearing in the proof of \Cref{pac_shell} by Pachner. (This is \Cref{collar}.)
\end{asparaenum}
\medskip
  Step 2 and 3 provide the proof of \Cref{analogIKNBoundary} by adapting the proof of Theorem 1.2 in \cite{IKN} to our setting.  
	
\subsection{Step 1: Restricting to manifolds with the same boundary}\label{sect:same boundary}
We consider two balanced PL homeomorphic manifolds $\Delta$ and $\Gamma$ of dimension $d$. Our aim is to show that, using shellings and balanced inverse shellings, we can transform them in such a way that they have isomorphic boundary complexes that moreover have the same induced coloring. 

First note that the boundary complexes $\partial \Delta$ and $\partial \Gamma$ are closed $(d-1)$-dimensional manifolds that are properly $(d+1)$-colorable. In fact, those boundaries might even be $d$-colorable and as such balanced. By \Cref{ikn_bist}, we know that there is a sequence of bistellar flips connecting $\partial\Delta$ with $\partial \Gamma$ such that each intermediate complex is properly $(d+1)$-colored. It now remains to encode this sequence of bistellar flips on $\partial\Delta$ as a sequence of shellings and balanced inverse shellings on $\Delta$. The next lemma, which we assume to be well-known, fulfills this task in the non-balanced situation. 
\begin{lemma}
	\label{bist_to_shell}
	Let $\Delta$ be a combinatorial $d$-manifold with boundary. Let $A\in \partial\Delta$ such that $\lk_{\partial\Delta}(A)=\partial B$ for some $B\notin \partial\Delta$. Then there exists a combinatorial $d$-manifold with boundary $\Delta'$ that is obtained from $\Delta$ by a single shelling, or inverse shelling, and $\partial\Delta'=\chi_{A,B}(\partial\Delta)$.
\end{lemma}

Combining the previous result with \Cref{ikn_bist}, we can now prove the main result of this section.

\begin{figure}[h]
	\centering
	\includegraphics[scale=0.8]{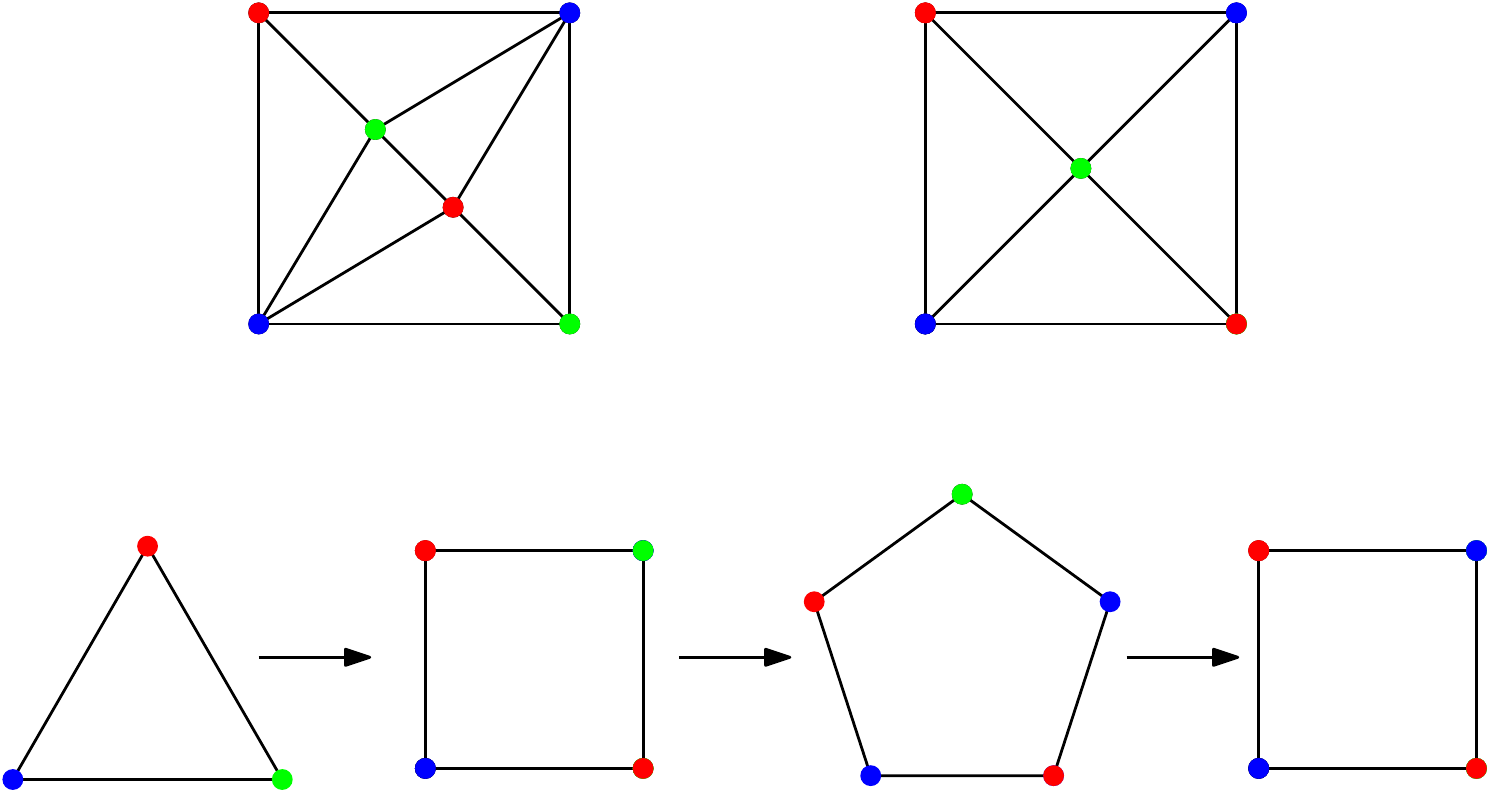}
	\caption{\emph{First row}: two balanced $2$-balls with isomorphic, though differently colored boundaries. \emph{Second row}: a sequence af bistellar flips connecting a triangle and a square with a prescribed coloring. }
	\label{same_bound_diff_color}
\end{figure}

\begin{proposition}
	\label{bound_same}
		Let $\Delta$ and $\Gamma$ be balanced combinatorial $d$-manifolds with boundary that are PL homeomorphic. Then there exists a balanced combinatorial $d$-manifold with boundary $\Delta'$ such that
		\begin{enumerate}
			\item $\Delta\bsh\Delta'$, 
			\item there exists a simplicial isomorphism $\varphi:\;\partial\Delta'\to\partial\Gamma$ that preserves the coloring, i.e., if $\kappa'$ and $\kappa$ are proper $(d+1)$-colorings of $\Delta'$ and $\Gamma$, respectively, then $\kappa'(v)=\kappa(\varphi(v))$ for all $v\in V(\partial(\Delta'))$.
		\end{enumerate}
\end{proposition}
\begin{proof}
	Since $\Delta$ and $\Gamma$ are balanced PL homeomorphic combinatorial $d$-manifolds with boundary, it follows that their boundaries, $\partial\Delta$ and $\partial\Gamma$, are closed combinatorial $(d-1)$-manifolds that are PL homeomorphic that are properly $(d+1)$-colorable. By \Cref{ikn_bist} there exists a sequence of bistellar flips from $\partial\Delta$ to $\partial\Gamma$ such that each intermediate complex is properly $(d+1)$-colored, and the flips preserve the vertex colors. Due to \Cref{bist_to_shell}, this sequence induces a sequence of shellings and inverse shellings from $\Delta$ to some PL homeomorphic manifold $\Delta'$ whose boundary is isomorphic to $\partial\Gamma$. Moreover, \Cref{ikn_bist} ensures that the vertex colors of $\partial\Delta'$ are preserved under this isomorphism. To see that none of the inverse shellings in the constructed sequence destroys balancedness, it is enough to remark that newly created edges lie in the boundary, which is itself properly $(d+1)$-colored. 
\end{proof}
In general, even if two balanced manifolds have the same boundary, their colorings restricted to the boundary might be different. The first row of \Cref{same_bound_diff_color} shows an example for this phenomenon. However, \Cref{ikn_bist} guarantees that one can enforce a particular coloring just using bistellar flips. An illustration for this is given in the second row of \Cref{same_bound_diff_color}, where the boundary of a triangle is related to the boundary of a square with a prescribed proper $2$-coloring through a sequence of bistellar flips. If we drop any requirement on the coloring, the first move already suffices.  
\Cref{bound_same} enables us to convert balanced manifolds $\Delta$ and $\Gamma$ into manifolds with the same boundary and then to glue those two manifolds along their boundary via the map $\varphi$. In this way, we will obtain a simplicial poset. This idea will be made more precise in the next section.

\subsection{Step 2: Constructing a shellable pseudo-cobordism}\label{Sect:pseudocobordism}
The aim of this section is to prove that analogs of Theorem 4.8, Lemma 5.2 and Corollary 5.3 of \cite{IKN} hold for combinatorial manifolds that have the same boundary. This will require to generalize the notion of a pseudo-cobordism from \cite{IKN}, that is for closed manifolds $\Delta$ and $\Gamma$, to non-closed manifolds with isomorphic boundaries. The proofs are almost verbatim the same as the ones in \cite{IKN} and we will therefore only describe the overall strategy and indicate differences.

\begin{figure}[h]
	\centering
	\includegraphics[scale=0.8]{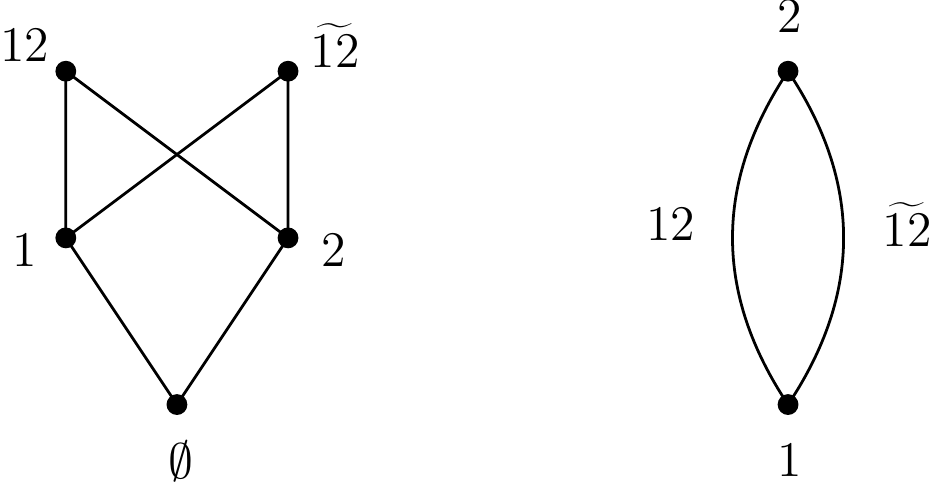}
	\caption{A simplicial poset that is not a simplicial complex (\emph{left}) and its corresponding geometric realization (\emph{right}).}
	\label{simp_poset}
\end{figure}

Before giving the definition of a pseudo-cobordism, we need to recall some notions concerning simplicial posets.  A \emph{simplicial poset} is a finite poset $\Omega$ with a unique minimal element $\emptyset$ such that for each $F\in\Omega$ the interval $\left[ \emptyset,F\right]=\{G\in \Omega~:~\emptyset \leq G\leq F\}$ is isomorphic to a Boolean lattice. Here, we denote with $\leq$ the order relation on $\Omega$.  A \emph{relative simplicial poset} is a pair of posets $(\Omega,\Sigma)$, such that $\Sigma\subseteq \Omega$ is a lower order ideal of $\Omega$, i.e., $\sigma\in \Sigma$ and $\tau<\sigma$ implies $\tau\in \Sigma$. Note that a simplicial poset $\Omega$ can be identified with the relative simplicial poset $(\Omega,\emptyset)$. The set of (relative) simplicial posets contains the set of (relative) simplicial complexes, but this inclusion is strict, see e.g., \Cref{simp_poset} for an example of a simplicial poset that is not a simplicial complex.

\emph{Faces} and \emph{facets} of a relative simplicial poset $(\Omega,\Gamma)$ are elements and inclusion-maximal elements of $\Omega\setminus \Gamma$, respectively. The \emph{dimension} of a face $F\in \Omega\setminus \Gamma$ is defined as $\dim F=\mathrm{rk}([\emptyset,F])-1$, where $\mathrm{rk}([\emptyset,F])$ denotes the rank of the Boolean interval $[\emptyset,F]$. The \emph{dimension} of $(\Omega,\Sigma)$ is the maximal dimension of its facets and we say that $(\Omega,\Sigma)$ is \emph{pure} if all facets have the same dimension. A pure relative simplicial poset $(\Omega,\Sigma)$ is called \emph{shellable} if there exists an ordering $F_1,\ldots,F_m$ of the facets of $(\Omega,\Sigma)$ such that for every $1\leq i\leq m$ the set
$$
\{G\in[\emptyset, F_i]~:~G\notin \bigcup_{j=1}^{i-1}[\emptyset,F_j]\cup\Sigma\}
$$
has a unique minimal element $R(F_i)$. 

The definition of a \emph{balanced} simplicial poset is completely analogous to the one of a balanced simplicial complex and is natural in the sense that a simplicial poset might have multiple edges but no loops. 
A $(d+1)$-dimensional simplicial poset $\Omega$ is a \emph{nonpure pseudomanifold} if every $d$-face is properly contained at most two facets. The \emph{pseudoboundary} $\tilde{\partial}\Omega$ of $\Omega$ is the subposet of $\Omega$ induced by the $d$-faces contained in zero or one $(d+1)$-face. We are now ready to define a pseudo-cobordism. 
\begin{definition}\label{def:pseudo-cobordism}
	Let $\Delta$ and $\Gamma$ be combinatorial $d$-manifolds with $\partial\Delta\cong\partial\Gamma$. (We also allow $\partial\Delta=\partial\Gamma=\emptyset$.) A \emph{pseudo-cobordism} $(\Omega,\varphi,\psi)$ between $\Delta$ and $\Gamma$ is a nonpure pseudomanifold $\Omega$ together with two simplicial embeddings $\varphi:\;\Delta\hookrightarrow\Omega$ and  $\psi:\;\Gamma\hookrightarrow\Omega$ such that:
	\begin{enumerate}
		\item[(1)] $\varphi(\Delta)\cup\psi(\Gamma)=\tilde{\partial}\Omega$,
		\item[(2)] a $d$-face $F\in\Omega$ lies in $\varphi(\Delta)\cap\psi(\Gamma)$ if and only if $F$ is not contained in any $(d+1)$-face of $\Omega$,
		\item[(3)] $\varphi(\partial\Delta)=\psi(\partial\Gamma)$.
	\end{enumerate}
\end{definition}
Note that, if $\Delta$ and $\Gamma$ are closed manifolds, then the assumption $\partial\Delta\cong\partial\Gamma$ is trivially satisfied and condition (3) is vacuous. In this case, we thus recover the definition of a pseudo-cobordism from \cite{IKN}. Also observe that condition (3) already implies $\partial\Delta\cong\partial\Gamma$, so that one might omit this assumption in the definition. However, we decided to keep it in order to emphasize that the definition is only for manifolds satisfying this condition. 
We will mostly be interested in \emph{shellable} pseudo-cobordisms:
\begin{definition}
Let $\Delta$ and $\Gamma$ be combinatorial $d$-manifolds such that $\partial\Delta\cong\partial\Gamma$. 
A pseudo-cobordism $(\Omega,\varphi,\psi)$ between $\Delta$ and $\Gamma$ is \emph{shellable} if there is an ordering $F_1,\dots,F_t$ of the $(d+1)$-faces of $\Omega$ such that
\begin{enumerate}
	\item $F_1,\dots,F_t$ is a shelling order on the relative simplicial poset $(\Omega,\varphi(\Delta))$,
	\item $F_t,\dots,F_1$ is a shelling order on the relative simplicial poset $(\Omega,\psi(\Gamma))$.
\end{enumerate}
\end{definition}
The simplest example of a shellable pseudo-cobordism is provided by a bistellar flip. More precisely, given a simplicial complex $\Delta$ and a face $A\in\Delta$ such that $\lk_{\Delta}(A)=\partial B$, for some $B\notin \Delta$, the simplicial complex $\Delta\cup(\langle A\rangle*\langle B\rangle)$ is called an \emph{elementary pseudo-cobordism}. Indeed, using that a bistellar flip does not modify the boundary, it is not difficult to see that if $\Delta$ is a combinatorial manifold, then $\Delta\cup(\langle A\rangle*\langle B\rangle)$ is a shellable pseudo-cobordism between $\Delta$ and $\chi_{A,B}(\Delta)$ with the obvious embeddings. The following characterization of a shellable pseudo-cobordism was shown in \cite[Proposition 4.7]{IKN}.

\begin{proposition}\label{prop:pseudocobordism composed}
A pseudo-cobordism is shellable if and only if it can be represented as a composition of elementary pseudo-cobordisms. 
\end{proposition}
 Though the statement in \cite{IKN} is only for pseudo-cobordisms between closed manifolds, their proof carries over verbatim to the situation of non-closed manifolds with isomorphic boundaries. Indeed, the ``Only-if''-part relies on a series of lemmas that only use part (1) and (2) of \Cref{def:pseudo-cobordism} but nowhere that the manifolds are assumed to be closed. The ``If''-part follows from the fact that the composition of shellable pseudo-cobordisms is again a shellable pseudo-cobordism \cite[Lemma 4.6]{IKN}. To see that this statement is true in our setting, it is enough to note that condition (3) of \Cref{def:pseudo-cobordism} is preserved under composition. 

As a corollary of \Cref{prop:pseudocobordism composed} one now obtains that Theorem 4.8 of \cite{IKN} remains true in our setting (with the same proof).
\begin{theorem}\label{thm:4.8}
Let $\Delta$ and $\Gamma$ be combinatorial manifolds such that $\partial\Delta\cong \partial \Gamma$. Then $\Delta$ and $\Gamma$ are bistellar equivalent if and only if there exists a shellable pseudo-cobordism between $\Delta$ and $\Gamma$. 
\end{theorem}
\begin{example}\label{example:pseudo}
Let $\Delta$ be a $1$-dimensional simplicial complex, consisting of $3$ consecutive edges. Let $\Gamma$ be the simplicial complex obtained from $\Delta$ by first applying a bistellar flip to the middle edge and then performing the inverse move, as depicted in the first row of \Cref{pseudocob}. Note that we have $\Delta=\Gamma$. This sequence of bistellar flips can be encoded in a shellable pseudo-cobordism between $\Delta$ and $\Gamma$. This is shown in the second row of \Cref{pseudocob}. In the bottom right picture, the complex $\varphi(\Delta)\cap\psi(\Gamma)$ is depicted in green, while the blue and the red segment are respectively $\varphi(\Delta)\setminus\psi(\Gamma)$ and  $\psi(\Gamma)\setminus\varphi(\Delta)$. We refer to \cite[Section 4]{IKN} for the precise construction of the pseudo-cobordism.
 \begin{figure}[h]
	\centering
	\includegraphics[scale=0.8]{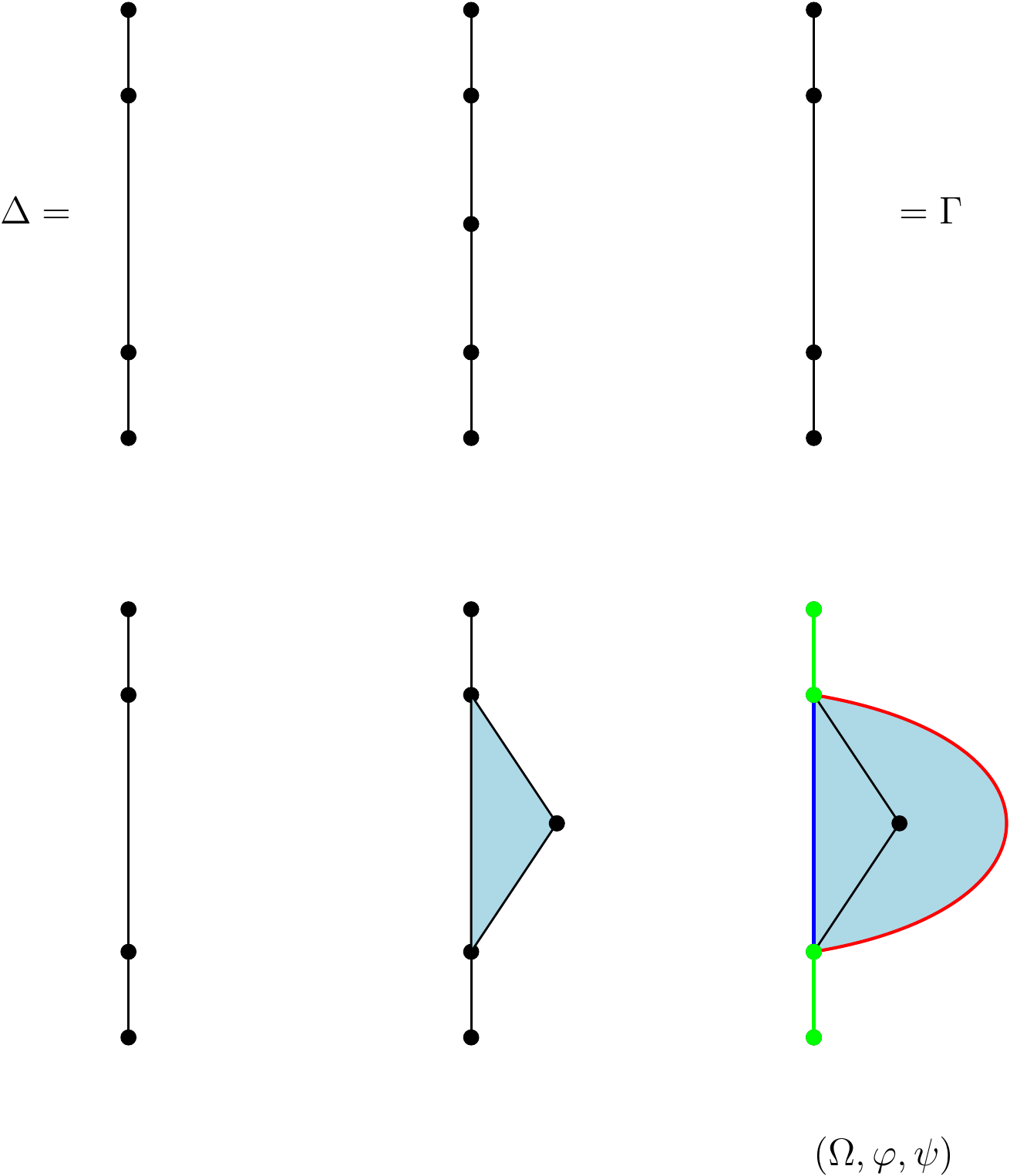}
	\caption{A sequence of 1-dimensional bistellar flips (\emph{first row}) encodes a sequence of elementary pseudo-cobordisms (\emph{second row}) and vice versa. }
	\label{pseudocob}
\end{figure}
\end{example}
In \cite[Corollary 5.3]{IKN} it is shown that, if there exists a shellable pseudo-cobordism between $\Delta$ and $\Gamma$, then it can always be chosen in such a way that $\Delta$ and $\Gamma$ embed disjointly. This is done by constructing a shellable pseudo-cobordism from $\Delta$ to some manifold $\Delta'$, such that $F\notin\Delta'$ for a given face $F\in\Delta$ (\cite[Lemma 5.2]{IKN}). Iterating this procedure over all vertices of $\Delta$, and then composing the obtained shellable pseudo-cobordism from $\Delta$ to $\Delta'$ with the shellable pseudo-cobordism between $\Delta'$ and $\Gamma$ yields the desired result. In our setting, we have the following analog: 
\begin{lemma}
Let $\Delta$ be a combinatorial $d$-manifold with boundary and let $F\in\overset{\circ}{\Delta}$ be a face in the interior of $\Delta$. Then there exists a combinatorial $d$-manifold $\Delta'$ with boundary and a $(d+1)$-dimensional nonpure pseudomanifold $\Omega$ such that 
\begin{enumerate}
\item[(1)] $\partial\Delta'\cong\partial\Delta$,
\item[(2)] $\Delta$ and $\Delta'$ are PL homeomorphic,
\item[(3)] $\Omega$ is a shellable pseudo-cobordism between $\Delta$ and $\Delta'$,
\item[(4)] $F\notin \Delta'$.
\end{enumerate}
\end{lemma}
We only comment on some specific points of the proof, since is basically the same as the one of \cite[Lemma 5.2]{IKN}.  There, the first ingredient is that the link of any face $F$ of a closed combinatorial $d$-manifold is a combinatorial $(d-|F|)$-sphere. In our setting this is true, since the face $F$ lies in the interior of $\Delta$. In the next step, the proof proceeds by constructing $\Delta'$ and $\Omega$ explicitly. In our situation, one immediately sees that this construction satisfies $\partial\Delta\cong \partial\Delta'$ and correspondingly for the embeddings, which is why (1) and (3) are hold. Finally, instead of using \Cref{pac_bist} and \cite[Theorem 4.8]{IKN}, we use  \Cref{same_bound_bist_equi} and \Cref{thm:4.8} to conclude that (2) holds.\\

Using the same arguments as in the proof of \cite[Corollary 5.3]{IKN}, we obtain the main result of this section.

\begin{corollary}\label{pseudo_disj}
Let $\Delta$ and $\Gamma$ be balanced PL homeomorphic combinatorial $d$-manifolds with boundary. Assume that $\partial\Delta\cong\partial\Gamma$ and that this isomorphism preserves the coloring. Then there exists a shellable pseudo-cobordism $(\Omega,\varphi,\psi)$ between $\Delta$ and $\Gamma$, such that $\varphi(\Delta)\cap\psi(\Gamma)=\varphi(\partial\Delta)=\psi(\partial\Gamma)$. Moreover, the pseudo-boundary $\tilde{\partial}\Omega$  is a balanced simplicial poset.
\end{corollary}
 The ``Moreover''-part is immediate from the fact that $\Delta$ and $\Gamma$ are both balanced and that their colorings coincide on their boundaries. We also want to remark that the pseudoboundary $\tilde{\partial}\Omega=\varphi(\Delta)\cup\psi(\Gamma)$ is a simplicial poset but not necessarily a simplicial complex. 

\subsection{Step 3: Converting bistellar flips into cross-flips}
The aim of this section is to prove \Cref{analogIKNBoundary}, providing an analog of Theorem 1.2 of \cite{IKN} for PL homeomorphic manifolds with isomorphic boundaries (see also \Cref{analogIKNBoundary}).
\begin{theorem}
	\label{cross_int}
	Let $\Delta$ and $\Gamma$ be balanced combinatorial $d$-manifolds with boundary that are PL homeomorphic. Assume moreover that $\partial\Delta\cong\partial\Gamma$ and that the colorings coincide on the boundary. Then, there exists a sequence of basic cross-flips connecting $\Delta$ and $\Gamma$.
\end{theorem}
The proof is analogous to the last steps of the one of Theorem 1.2 in \cite{IKN}. We include it as a service to the reader. 
\begin{proof}
	We apply \Cref{pseudo_disj} to obtain a shellable pseudo-cobordism $(\Omega,\varphi,\psi)$ between $\Delta$ and $\Gamma$ with $\varphi(\Delta)\cap\psi(\Gamma)=\varphi(\partial\Delta)=\psi(\partial\Gamma)$. Though $\Omega$ might not be balanced, it follows from \cite[Corollary 3.2]{IKN} that there exists a balanced nonpure $(d+1)$-pseudomanifold $\Omega'$ obtained by stellar subdivision of interior faces of $\Omega$, which is a pseudo-cobordism between $\Delta$ and $\Gamma$. Since stellar subdivisions preserve shellability, $\Omega'$ is shellable. (This follows from Proposition 5.7 in \cite{IKN}, which is only stated for closed manifolds but whose proof also goes through in our setting.) Applying the diamond operation to the simplicial poset $\Omega'$ leads to a cross-polytopal $(d+1)$-complex $\Diamond(\Omega')$, and the shelling order on $(\Omega',\varphi,\psi)$ induces an order on the $(d+1)$-cells of $\Diamond(\Omega')$, which encodes a sequence of cross-flips between $\Delta$ and $\Gamma$. For more details on this part, we refer to the proof of Theorem 3.10 in \cite{IKN} and to the next example as an illustration.
\end{proof}	

\begin{example}
As in \Cref{example:pseudo} we let $\Delta=\Gamma$ be the $1$-dimensional ball consisting of $3$ consecutive edges whose vertices are colored alternately with red and blue. The pseudo-cobordism provided in \Cref{example:pseudo} does not satisfy the conditions in \Cref{pseudo_disj} since the embeddings of $\Delta$ and $\Gamma$ do not intersect just in two vertices but two edges. However using stellar subdivisions we can construct a shellable pseudo-cobordism $(\Omega,\varphi,\psi)$ that meets the requirements of \Cref{pseudo_disj}. Such a pseudo-cobordism, though not a minimal one, is depicted in the top left image of \Cref{pro}. The labels on the $2$-faces encode a shelling order $F_1,\dots,F_{10}$ of $(\Omega,\varphi(\Delta))$. Applying the diamond operation, we obtain a cross-polytopal complex. In this case, the diamond operation subdivides all $1$-faces not containing vertices of color 0 (which is blue in the figure), i.e., all edges whose vertices are colored with red and green. The $2$-cells of $\Diamond(\Omega)$ are $2$-dimensional cross-polytopes. We now describe how an inverse shelling gets encoded into a cross-flip, taking $F_3$ as an example. If we identify $\Diamond(F_3)$ with $\mathcal{C}_2$ and consider the decomposition $F_3=A_3\cup R_3$ given by the restriction face $R_3$ and $A_3=F_3\setminus R_3$, we notice that adding $F_3$ to the simplicial poset corresponds to replacing $\Diamond(A_3)$ with $\Diamond(\langle R_3\rangle *\partial A_3)=\mathcal{C}_2\setminus\Diamond(A_3)$, which describes the application of the cross-flip $(\Diamond(A_3),\C_2\setminus\Diamond(A_3))$, see the bottom right of \Cref{pro}).
\begin{figure}[h]
	\centering
	\includegraphics[scale=0.8]{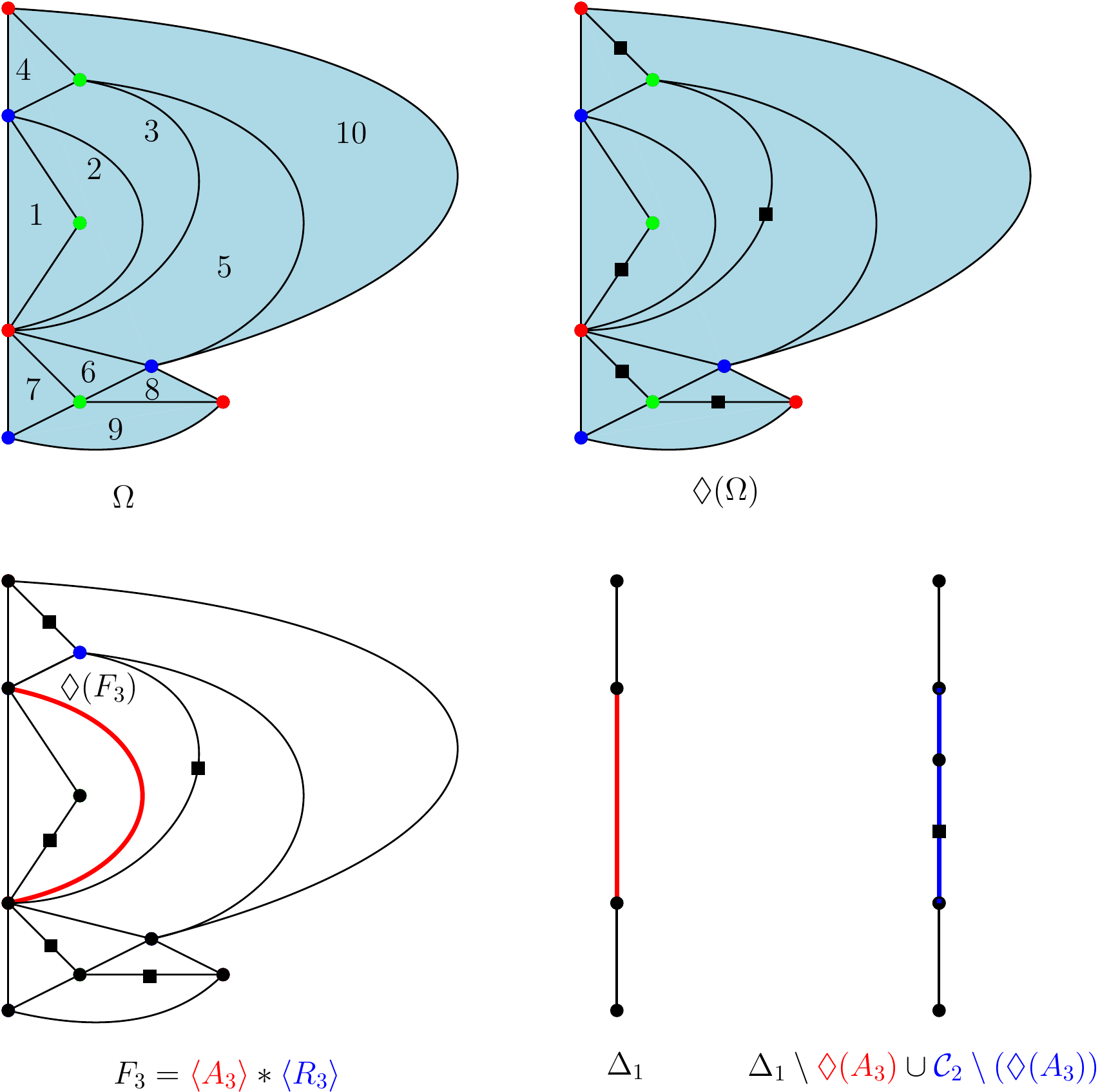}
	\caption{\emph{First row}: a shellable pseudo-cobordism $(\Omega,\varphi,\psi)$ and the cross-polytopal complex $\Diamond(\Gamma)$. \emph{Second row}: A shelling on $\Omega$ encodes a cross-flip between the intermediate steps. }
	\label{pro}
\end{figure}
\end{example}

\subsection{Step 4: Building a collar and shelling cross-flips}\label{sect:collar and shelling}
The last ingredient, we need for the proof of \Cref{mainresult}, is a way to convert every cross-flip into a sequence of shellings and balanced inverse shellings. Having achieved this, the basic strategy is the following: Whenever a cross-flip is performed in the interior of a manifold, we will first ``enter'' it from an arbitrary boundary facet, ``dig'' into the manifold by removing facets (using shellings) until we meet a facet that is involved in the cross-flip to be carried out. In the next step, we shell the subcomplex to be removed, add its complement with respect to the boundary of the cross-polytope using inverse shellings and finally close the path, we built before, using inverse shellings. To avoid weird and undesirable side effects, when shelling the cross-flip, we need to make sure that the only facet involved in the cross-flip that meets the boundary of the manifold in codimension 1 is the facet that was first hit when digging the tunnel into the manifold. This can be achieved by building a collar around the manifold (using shellings and balanced inverse shelling), an idea going back to the proof of \Cref{pac_shell} by Pachner. 

The following results will be crucial.
 
\begin{lemma}\cite[Theorem 5.8]{Pac1}
	\label{shell_ball}
	Every combinatorial $d$-sphere $\Delta$ is the boundary of a shellable combinatorial $(d+1)$-ball $\Omega$. Moreover, $\Omega$ can be chosen so that $\Delta$ is an induced subcomplex of $\Omega$.
\end{lemma}

\begin{theorem}\cite[Theorem 3.1]{IKN}
	\label{relative_ikn}
	Let $\Omega$ be a $d$-dimensional simplicial complex and $\Delta\subsetneq \Omega$ a subcomplex. Let $\kappa:V(\Delta)\longrightarrow \left\lbrace 0,\dots,m-1\right\rbrace$ be a proper $m$-coloring of $\Delta$.Then there is a stellar subdivision $\Omega'$ of $\Omega$ s.t.:
	\begin{enumerate}
		\item[(1)] $\Delta$ is a subcomplex of $\Omega'$.
		\item[(2)] The coloring $\kappa$ extends to a proper coloring $\kappa':V(\Omega')\longrightarrow\left\lbrace 0,\dots, \max\left\lbrace m-1,d\right\rbrace\right\rbrace $ such that all vertices not in $\Delta$ receive colors in $\left\lbrace 0,\dots,d\right\rbrace $.
	\end{enumerate} 
\end{theorem}

\begin{remark}
	\label{special_ball}
 Let $\Delta$ be a combinatorial $d$-sphere and let $\kappa:V(\Delta)\to \{0,\ldots,m-1\}$ be a proper $m$-coloring of $\Delta$. Then, by \Cref{shell_ball}, there exists a shellable $(d+1)$-ball $\Omega$ with $\partial\Omega=\Delta$ and such that $\Delta$ is an induced subcomplex of $\Omega$. \Cref{relative_ikn} further yields a $(d+1)$-ball $\Omega'$ such that:
	\begin{enumerate}
		\item $\partial\Omega'=\Delta$,
		\item $\Omega'$ is shellable,
		\item $\Omega'$ is properly $(\max\left\lbrace m,d+2\right\rbrace)$-colored and this extends the coloring $\kappa$ of $\Delta$,
		\item $\Delta$ is an induced subcomplex of $\Omega'$.
	\end{enumerate}
	Conditions (2) and (4) hold since stellar subdivisions preserve both shellability and the property of being an induced subcomplex. 
\end{remark}

We state a further lemma.
\begin{lemma}\cite[Lemma 4.7]{Pac1}
	\label{pachner_link}
	Let $\Delta$ be a combinatorial $d$-manifold and $\Omega\subseteq\partial\Delta$ a shellable $(d-1)$-ball such that $\Omega\subseteq\lk_{\Delta}(v)$ for some vertex $v\in\overset{\circ}{\Delta}$. Then there is a sequence of shellings converting $\Delta$  into $\Delta\setminus( \langle v\rangle*\overset{\circ}{\Omega})$.
\end{lemma}
Putting the previous results together, allows us to build a collar:  

\begin{theorem}
	\label{collar}
	Let $\Delta$ be a balanced combinatorial $d$-manifold with boundary and let $F$ be a facet of $\partial \Delta$. Then there exists a balanced combinatorial $d$-manifold $\Delta'$ such that
	\begin{enumerate}
		\item $\Delta'$ can be transformed into $\Delta$ by a sequence of shellings,
		\item $\Delta$ is an induced subcomplex of $\Delta'$,
		\item $\partial \Delta\cap \partial \Delta'=\left\langle F\right\rangle $.
	\end{enumerate}
\end{theorem}
\begin{proof}
We fix a coloring $\kappa$ of $\Delta$. 
Let $v\in V(\partial\Delta)$ with $v\notin F$. Then $\lk_{\partial\Delta}(v)$ is a combinatorial $( d-2) $-sphere, which is properly $d$-colored (though not necessarily balanced). By \Cref{special_ball} we can choose a balanced, shellable $( d-1) $-ball $\Omega$, whose boundary is $\lk_{\partial \Delta}(v)$ and the latter is an induced subcomplex of $\Omega$. Moreover, there is a proper $d$-coloring $\kappa'$ of $\Omega$ that restricts to $\kappa$ on $\lk_{\partial\Delta}(v)$. Let us consider $\widetilde{\Delta}:=\Delta\cup\left( \langle v\rangle *\Omega\right)$. Since $\partial\widetilde{\Delta}= (\partial \Delta\setminus\left( \st_{\partial\Delta}(v)\right)) \cup\Omega$, the open star $\overset{\circ}{\st_{\partial\Delta}(v)}$ is contained in $\overset{\circ}{\widetilde{\Delta}}$. As $\lk_{\partial\Delta}(v)$ is an induced subcomplex of $\Omega$, the same is true for $\Delta$, considered as a subcomplex of $\widetilde{\Delta}$. Since the coloring $\kappa'$ of $\Omega$ does not use $\kappa(v)$ and since $\kappa'|_{\lk_{\partial\Delta}(v)}=\kappa$, we conclude that $\Omega\ast \langle v\rangle$ and thus also $\widetilde{\Delta}$ is balanced. Finally, it follows from \Cref{pachner_link}, that there is a sequence of shellings from $\widetilde{\Delta}$ to $\Delta$. We now apply the described construction to every vertex $v\in V(\partial\Delta)\setminus F$ to obtain a simplicial complex $\Delta'$. As any face $G\in \partial\Delta$ with $G\not\subseteq F$ lies in the open star $\overset{\circ}{\st_{\partial\Delta}(v)}$ of some vertex $v\in V(\partial\Delta)\setminus F$, it follows that $\Delta'$ satisfies condition (3).
\end{proof}

As the only missing ingredient for the proof of \Cref{mainresult} we need to convert cross-flips into shellings and their inverses in such a way that balancedness is preserved. This is done in the remaining part of this section.\\

\begin{figure}[h]  
	\centering
	\includegraphics[scale=0.7]{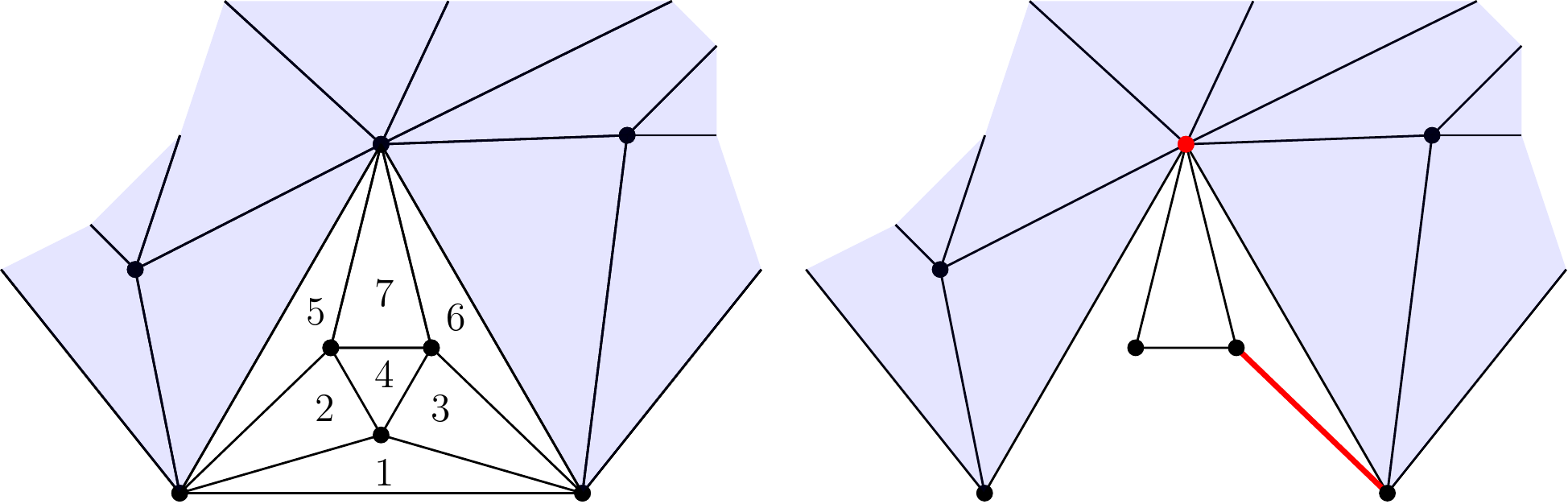}
	\caption{A shelling order on $D$ ({\em left}) that is not a of shelling for $(\Delta,\Delta\setminus D)$ ({\em right}).}	\label{non-shelling}
\end{figure}

In the following, we assume that $\Delta$ is a combinatorial $d$-manifold with boundary, $D=\Diamond(\Gamma)\subseteq \Delta$ is an induced subcomplex for some $d$-ball $\Gamma\subseteq \partial\sigma^{d+1}$ and $D\cap \partial\Delta=\langle F\rangle$ for a $(d-1)$-face $F\subseteq \partial\Delta$. We let $\Delta'=\chi^*_D(\Delta)$. Our aim is to construct a sequence of shellings from $\Delta$ to $\Delta\setminus D$ and a sequence of inverse shellings from $\Delta\setminus D$ to $\Delta'=\Delta\setminus D\cup (\C_d\setminus D)$ that preserves balancedness. Since we know that $D$ and $\C_d\setminus D$ are both shellable, an obvious choice might be a reversed shelling order on $D$, followed by a shelling order on $\C_d\setminus D$. However, in general this will not work, since faces of $D$ and $\C_d\setminus D$ might also intersect $\Delta\setminus D$, which can cause obstructions to shellability. An instance for this behavior is given in \Cref{non-shelling}. On the one hand, the ordering of the facets, indicated by their labels, is a reversed shelling for the designated subcomplex. On the other hand, taking into account the large complex, we are not allowed to remove the facet labelled $6$ (once facets $1,\ldots, 5$ have been removed) since it intersects the boundary of the given manifold in a nonpure subcomplex, consisting of an edge and an isolated vertex (shown in red in the right picture of \Cref{non-shelling}). So, what we need to construct is a shelling of the relative simplicial complex $(\Delta,\Delta\setminus D)$, which is defined as for relative simplicial posets (see \Cref{Sect:pseudocobordism}).

In the sequel, we will define an ordering on the facets of $D$ and $\C_d\setminus D$, that provides a sequence of shellings for $\Delta$ and $\Delta'$, respectively. In this way, we can relate $\Delta$ and $\Delta'$  by the sequence for $D$ followed by the reversed sequence for $\C_d\setminus D$. First, we need some preparations. Let $\sigma^{d+1}$ be the $(d+1)$-simplex on vertex set $\{0,1,\ldots,d+1\}$ and for $0\leq i\leq d+1$ let $\Gamma_i=\{0,\ldots,\hat{i},\ldots,d+1\}$ be the facet not containing the vertex $i$. We start with a description of the facets of $\Diamond(\Gamma_i)$ for $0\leq i\leq d$, where we follow the notation of the paragraph preceding \Cref{def:basic cross-flips}. 
\begin{lemma}\label{lemma:diamond complexes}
For $0\leq i\leq d-1$ let $\C_{d-i-1}(i+1,\ldots,d)$ be the boundary complex of the $(d-i-1)$-dimensional cross-polytope on vertex set $\{i+1,\ldots,d\}\cup\{v_{i+1},\ldots,v_d\}$ and set $\C_{-1}=\{\emptyset\}$. Then,
\begin{equation*}
\Diamond(\Gamma_i)=\begin{cases}
\left\langle \{0,\ldots,i-1,v_i\}\right\rangle \ast\C_{d-i-1}(i+1,\ldots,d) &\mbox{ for }  0\leq i\leq d,\\
\left\langle \{0,\ldots,d\}\right\rangle  &\mbox{ for } i=d+1.
\end{cases}
\end{equation*}
\end{lemma}

\begin{figure}[h]
	\centering
	\includegraphics[scale=0.8]{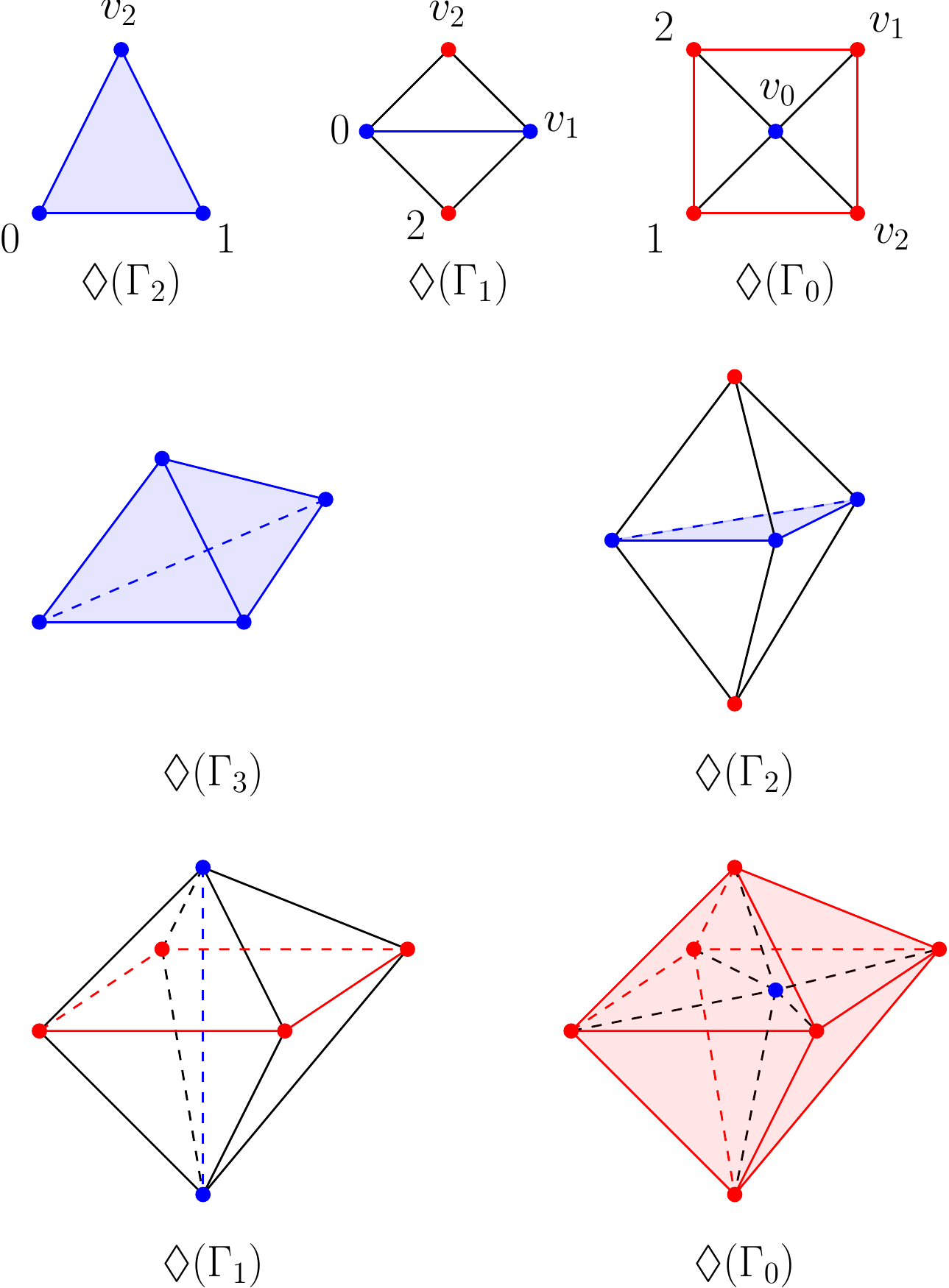}
	\caption{All the subcomplexes $\Diamond(\Gamma_i)$ for $d=2$ (\emph{first row}) and $d=3$ (\emph{second row}). The coloring here indicates the decomposition $\Diamond(\Gamma_i)={\color{blue}\sigma^{i}}*{\color{red}\mathcal{C}_{d-i-1}(i+1,\ldots,d)}$. }
	\label{diamonds}
\end{figure}

\begin{proof}
First note that the facet $\Gamma_{d+1}$ is not subdivided by the diamond operation.\\
Let $0\leq i\leq d$. In this case, the diamond operation will successively subdivide the faces $F_j=\{j+1,\ldots,d+1\}$ for $i\leq j\leq d$, i.e.,
$$
\Diamond(\Gamma_i)=\mathrm{sd}_{F_{d}}\circ\mathrm{sd}_{F_{d-1}}\circ\cdots \circ\mathrm{sd}_{F_{i}}(\Gamma_i).
$$
We prove the more general statement that for any $i\leq j\leq d$,
\begin{align*}
\mathrm{sd}_{F_{j}}\circ\mathrm{sd}_{F_{j-1}}\circ\cdots \circ\mathrm{sd}_{F_{i}}(\Gamma_i)=\left\langle \{0,\ldots,i-1,v_i\}\right\rangle \ast\partial\left\langle \{v_{i+1},i+1\}\right\rangle \ast \cdots\ast\partial\left\langle \{v_j,j\}\right\rangle \ast\partial F_j
\end{align*}
We proceed by induction on $j$. For $j=i$, we have
$$
\mathrm{sd}_{F_i}(\Gamma_i)=\left\langle \{0,\ldots,i-1\}\right\rangle \ast \left\langle \{v_i\}\right\rangle \ast \partial F_{i},$$
which implies the desired statement. For $j>i$, it holds that 
\begin{align*}
&\mathrm{sd}_{F_{j}}\circ\mathrm{sd}_{F_{j-1}}\circ\cdots \circ\mathrm{sd}_{F_{i}}(\Gamma_i)\\
=&\mathrm{sd}_{F_j}(\left\langle \{0,\ldots,i-1,v_i\}\right\rangle \ast\partial\left\langle \{v_{i+1},i+1\}\right\rangle \ast \cdots\ast\partial\left\langle \{v_{j-1},j-1\}\right\rangle \ast \partial F_{j-1})\\
=&\left\langle \{0,\ldots,i-1,v_i\}\right\rangle \ast\partial\left\langle \{v_{i+1},i+1\}\right\rangle \ast \cdots\ast\partial\left\langle \{v_{j-1},j-1\}\right\rangle \ast\mathrm{sd}_{F_j}(\partial F_{j-1})\\
=&\left\langle \{0,\ldots,i-1,v_i\}\right\rangle \ast\partial\left\langle \{v_{i+1},i+1\}\right\rangle \ast \cdots\ast\partial\left\langle \{v_{j-1},j-1\}\right\rangle \ast\partial\left\langle \{v_j,j\}\right\rangle \ast\partial F_j,
\end{align*}
where we use that 
$$
\mathrm{sd}_G(\Gamma\ast \Delta)=\Gamma\ast\mathrm{sd}_G(\Delta) \mbox{ for } G\in \Delta
$$
and
$$
\mathrm{sd}_{F_j}(\partial F_{j-1})=\partial\left\langle \{v_j,j\}\right\rangle \ast\partial F_j.
$$
The claim now follows.
\end{proof}
\Cref{diamonds} shows an illustration of the decomposition of $\Diamond(\Gamma_i)$ as provided by \Cref{lemma:diamond complexes}.

\begin{remark}
As a consequence of \Cref{lemma:diamond complexes} we can describe the boundary complex of $\Diamond(\Gamma_i)$ as 
$$
\partial\Diamond(\Gamma_i)=\begin{cases}
\partial\left\langle \{0,\ldots,i-1,v_i\}\right\rangle \ast\C_{d-i-1}(i+1,\ldots,d) &\mbox{ for }  0\leq i\leq d,\\
\partial\left\langle \{0,\ldots,d\}\right\rangle  &\mbox{ for } i=d+1.
\end{cases}
$$
Consequently, it follows that for $0\leq i<k\leq d+1$,
\begin{align*}
&\partial\Diamond(\Gamma_i)\cap \partial\Diamond(\Gamma_k)\\
=
&\begin{cases}
\left\langle \{0,\ldots,i-1,i+1,\ldots,k-1,v_k\}\right\rangle \ast\C_{d-k+1}(k+1,\ldots,d) &\mbox{ for } 1\leq k\leq d,\\
\partial\left\langle \{0,\ldots,i-1,i+1,\ldots,d\}\right\rangle  &\mbox{ for } k=d+1.
 \end{cases}
\end{align*}
In particular, $\Diamond(\Gamma_i)$ and $\Diamond(\Gamma_k)$ intersect in a pure $(d-1)$-dimensional subcomplex of their boundaries. Since any facet $F$ in $\Diamond(\Gamma_k)$ is of the form $\{0,\ldots,k-1,v_k\}\cup F'$ for a facet $F'$ of $\C_{d-k-1}(k+1,\ldots,d)$, we can further conclude that for any $F\in\mathcal{F}(\Diamond(\Gamma_k))$ and any $0\leq i<k\leq d+1$,  there exists $G\in \mathcal{F}(\Diamond(\Gamma_i))$ such that $\dim(F\cap G)=d-1$. Namely, $G=F\setminus\{i\}\cup\{v_i\}$. Observe, that we can further infer that for any $0\leq k<i\leq d+1$, there exists a facet $F\in \Diamond(\Gamma_k)$ and a facet $G\in\Diamond(\Gamma_i)$ such that $\dim (F\cap G)=d-1$. 
\end{remark}

As before, let $\Delta$ be a balanced combinatorial $d$-manifold with boundary. Let $\Gamma=\langle \Gamma_{i_1},\ldots,\Gamma_{i_k}\rangle$, where the $i_j$ are pairwise distinct and $i_2<\cdots <i_k$, and let $D=\Diamond(\Gamma)$ be an induced subcomplex of $\Delta$ such that $\partial\Delta\cap D=\langle F\rangle$ for a $(d-1)$-face $F$. Without loss of generality, assume $F\in \Diamond(\Gamma_{i_1})$. Let $G$ be the unique facet of $\Diamond(\Gamma_{i_1})$ containing $F$. 
We now describe a shelling for $(\Delta,\Delta\setminus D)$, starting with the facets of $\Diamond(\Gamma_{i_1})$ followed by the facets of $\Diamond(\Gamma_{i_2})$, $\ldots$, $\Diamond(\Gamma_{i_k})$.\\
We need some further notation. For $2\leq \ell \leq k$, we let $1\leq m(\ell)\leq \ell$ such that $i_{m(\ell)}=\min\{i_j~:~1\leq j\leq \ell \mbox{ and } i_j\geq i_\ell\}$. 
For $1\leq \ell\leq k$, we define the \emph{initial facet} $F_{\mathrm{in}}^{(\ell)}$ of $\Diamond(\Gamma_{i_\ell})$ (with respect to $\Gamma$ and $F$) as 
$$
F_{\mathrm{in}}^{(\ell)}=\begin{cases}
G&\mbox{ for } \ell=1,\\
\{0,\ldots,i_\ell-1,v_{i_\ell}\}\cup\{i_{\ell}+1,\ldots,i_{m(\ell)}-1\}\cup\{v_{i_{m(\ell)}},\ldots,v_d\}&\mbox{ for } 2\leq \ell\leq k.
\end{cases}
$$

We have the following simple observation:
\begin{lemma}\label{lemma:free facet}
For any $1\leq \ell\leq k$, the initial facet $F_{\mathrm{in}}^{(\ell)}$ of $\Diamond(\Gamma_{i_\ell})$ intersects the boundary of $\Delta_{\ell}=\Delta\setminus \left(\bigcup_{j=1}^{\ell-1}\Diamond(\Gamma_{i_j})\right)$ in a pure subcomplex of dimension $d-1$. In particular, the operation $\Delta_\ell\mapsto \Delta_\ell\setminus F_{\mathrm{in}}^{(\ell)}$ is an elementary shelling.
\end{lemma}

\begin{proof}
The statement is clear for $i_1$.\\
Let $\ell \geq 2$. We have 
\begin{align*}
 &F_{\mathrm{in}}^{(\ell)}\cap \left(\Delta\setminus \left(\bigcup_{j=1}^{\ell-1}\Diamond(\Gamma_{i_j})\right)\right)=F_{\mathrm{in}}^{(\ell)}\cap \partial\left(\bigcup_{j=1}^{\ell-1}\Diamond(\Gamma_{i_j})\right)\\
=& F_{\mathrm{in}}^{(\ell)}\cap \left(\bigcup_{j=1}^{\ell-1}\partial\Diamond(\Gamma_{i_j})\right)
=\bigcup_{j=1}^{\ell-1}\left( F_{\mathrm{in}}^{(\ell)}\cap \partial\Diamond(\Gamma_{i_j})\right),
\end{align*}
where the second equality follows from the fact that any $(d-1)$-face of $\Delta$ is contained in at most 2 facets. 
First assume $m_\ell\neq \ell$ and hence $i_{m_\ell}>i_\ell$. In this case,
\begin{align*}
 F_{\mathrm{in}}^{(\ell)}\cap \partial\Diamond(\Gamma_{i_j})=\begin{cases}
F_{\mathrm{in}}^{(\ell)}\setminus \{i_j\} &\mbox{ if } i_j<i_\ell \tag{\theequation}\label{eq:Intersection}\\ 
F_{\mathrm{in}}^{(\ell)}\setminus \{v_{i_\ell}\} &\mbox{ if } i_j=i_{m(\ell)}\\
F_{\mathrm{in}}^{(\ell)}\setminus \{v_{i_\ell},v_{i_{m(\ell)}},\ldots,v_{i_j-1}\}  &\mbox{ if } i_j>i_{m(\ell)}.
\end{cases}
\end{align*}
As the intersection in the last case is clearly contained in $F_{\mathrm{in}}^{(\ell)}\cap \partial\Diamond(\Gamma_{i_{m(\ell)}})$, the claim follows. \\
If $m(\ell)=\ell$, then $i_j<i_\ell$ for all $1\leq j\leq \ell-1$ and the claim follows from \eqref{eq:Intersection}.\\ 
The ``In particular''-part is now immediate (see e.g., \cite[Definition 5.1.11]{BH-book}).
\end{proof}

\Cref{lemma:free facet} in particular implies that for any $1\leq \ell\leq k$ the initial facet $F_{\mathrm{in}}^{(\ell)}$ is \emph{free} with respect to $\Delta\setminus\left(\bigcup_{j=1}^{\ell-1}\Diamond(\Gamma_{i_j})\right)$, meaning that it intersects the boundary of $\Delta\setminus\left(\bigcup_{j=1}^{\ell-1}\Diamond(\Gamma_{i_j})\right)$ in a $(d-1)$-ball.

We can finally define an ordering on the facets of $\Diamond(\Gamma_{i_\ell})$ with respect to the initial facet $F_{\mathrm{in}}^{(\ell)}$. This ordering will be inspired by the lexicographic ordering, that provides a shelling for $\C_{d}$ (see e.g., \cite[Theorem 1.14]{Longueville}). By \Cref{lemma:diamond complexes}, any facet $F\in \Diamond(\Gamma_{i_\ell})$ is of the form 
$\{0,\ldots,i_\ell-1,v_{i_\ell}\}\cup F'$, where $F'$ is a facet of $\C_{d-i_{\ell}-1}(i_\ell+1,\ldots,d)$. We can therefore interpret $F$ as an ordered $(d+1)$-tuple in $\{0\}\times \cdots \times \{i_{\ell-1}\}\times\{v_{i_\ell}\}\times\{i_\ell+1,v_{i_\ell+1}\}\times \cdots \times \{d,v_d\}$. We write $F(j)$ for the $j$\textsuperscript{th} entry of $F$, i.e., $F(j)=F\cap\{j,v_j\}$ for $0\leq i\leq d$. We use the same notation for tuples in $\{0,1\}^{d-i_\ell}$. We define the \emph{characteristic function} 
\begin{align*}
 \varphi_F:\;\mathcal{F}\left(\Diamond(\Gamma_{i_\ell})\right)&\to \left\lbrace0,1 \right\rbrace^{d-i_\ell}
\end{align*}
on $\mathcal{F}(\Diamond(\Gamma_{i_\ell}))$ with respect to $F$ by setting
$$
\varphi_F(G)(j)=\begin{cases}
0\quad \mbox{ if } G(i_\ell+1+j)=F(i_\ell+1+j),\\
1\quad\mbox{ if } G(i_\ell+1+j)\neq F(i_\ell+1+j)
\end{cases}
$$
for $0\leq j\leq d-i_\ell-1$. We set $\deg_F(G):=\sum_{j=1}^{d_{i_\ell}}\varphi_F(G)(j)$ and call this the \emph{degree} of $G$ with respect to $F$. As $\deg_F(G)$ counts the number of elements in $G\setminus (F\cap G)$, we can interpret $\deg_F(G)$ as a measure of similarity between $F$ and $G$. 
We now order the facets of $\Diamond(\Gamma_{i_\ell})$ by ordering $\{\varphi_{F_{\mathrm{in}}^{(\ell)}}(F)~:~F\in \Diamond(\Gamma_{i_\ell})\}$ according to the degree lexicographic ordering, i.e., we set
 $F\prec F'$  if 
\begin{itemize}
\item[(i)] $\deg_{F_{\mathrm{in}}^{(\ell)}}(\varphi(F))<\deg_{F_{\mathrm{in}}^{(\ell)}}(\varphi(F'))$, or 
\item[(ii)] $\deg_{F_{\mathrm{in}}^{(\ell)}}(\varphi(F))=\deg_{F_{\mathrm{in}}^{(\ell)}}(\varphi(F'))$  and there exists $m$ such that $\varphi_{F_{\mathrm{in}}^{(\ell)}}(F)(k)=\varphi_{F_{\mathrm{in}}^{(\ell)}}(F')(k)$ for all $1\leq k\leq m$ and $\varphi_{F_{\mathrm{in}}^{(\ell)}}(F)(m+1)<\varphi_{F_{\mathrm{in}}^{(\ell)}}(F')(m+1)$.
\end{itemize}
 In the latter case, we must have $\varphi_{F_{\mathrm{in}}^{(\ell)}}(F)(m+1)=0$ and $\varphi_{F_{\mathrm{in}}^{(\ell)}}(F')(m+1)=1$. For simplicity, we call this ordering \emph{degree lexicographic ordering} with respect to $ F_{\mathrm{in}}^{(\ell)}$. 

\begin{theorem}\label{theorem:shelling}
Let $\Delta$ be a balanced combinatorial $d$-manifold with boundary. Let $\Gamma=\langle \Gamma_{i_1},\ldots,\Gamma_{i_k}\rangle$, where the $i_j$ are pairwise distinct and $i_2<i_3<\cdots<i_k$, and let $\Diamond(\Gamma)$ be an induced subcomplex of $\Delta$ intersecting $\partial \Delta$ in a $(d-1)$-face $F\in \Diamond(\Gamma_{i_1})$. For $1\leq \ell \leq k$, let $r_\ell=|\mathcal{F}(\Diamond(\Gamma_{i_\ell}))|-1$ and let 
$$
F_0^{(\ell)}=F_{\mathrm{in}}^{(\ell)},F_1^{(\ell)},\ldots,F_{r_\ell}^{(\ell)}
$$ 
be the ordering of the facets of $\Diamond(\Gamma_{i_\ell})$ according to the degree lexicographic ordering with respect to $F_{\mathrm{in}}^{(\ell)}$. Then
$$
F_{r_k}^{(k)},\ldots,F_1^{(k)},F_{0}^{(k)},F_{r_{k-1}}^{(k-1)},\ldots,F_1^{(k-1)},F_{0}^{(k-1)},\ldots F_{r_0}^{(1)},\ldots, F_1^{(1)},F_{0}^{(1)}
$$
is a shelling on $\left(\Delta,\Delta\setminus \left(\bigcup_{j=1}^{k}\Diamond(\Gamma_{i_j})\right)\right)$.
\end{theorem}

\begin{proof}
We show that for any $1\leq \ell \leq k$ the ordering $F_{r_\ell}^{(\ell)},\ldots,F_1^{(\ell)},F_{0}^{(\ell)}$ is a shelling for the relative complex $\left(\Delta\setminus \left(\bigcup_{j=1}^{\ell-1}\Diamond(\Gamma_{i_j})\right),\Delta\setminus \left(\bigcup_{j=1}^{\ell}\Diamond(\Gamma_{i_j})\right)\right)$. 
We treat the cases $\ell=1$ and $\ell>1$ separately, since they require slightly different arguments. However, a general remark is in order first. Let $\Delta'$  be a pure subcomplex of $\Delta$ with $\Delta\setminus \Diamond(\Gamma)\subseteq \Delta'$  and let $A\in \Diamond(\Gamma)\cap \Delta'$. Since $\partial\Delta\cap \Diamond(\Gamma)=\langle F\rangle$, the face $A$ is in the interior of $\Delta'$ if and only if it is only contained in facets of $\Delta$ that also belong to $\Delta'$.

First assume $\ell=1$. For $0\leq i\leq r_1$, we set $\Delta_i=\Delta\setminus \langle F_0^{(1)},\ldots, F_{i-1}^{(1)}\rangle$. We need to show that the operation $\Delta_{i}\mapsto \Delta_i\setminus F_i^{(1)}=\Delta_{i+1}$ is an elementary shelling. For $F_0^{(1)}$ this follows from \Cref{lemma:free facet}. For $1\leq i\leq r_1$, we set
\begin{equation*}
A_i^{(1)}=
\{F_i^{(1)}(j)~:~\varphi_{F_{\mathrm{in}}^{(1)}}(F_i^{(1)})(i_1+1+j)=1\}
\end{equation*}
and
\begin{equation*}
R_i^{(1)}=\{0,\ldots,i_1-1,v_{i_1}\}\cup\{F_i^{(1)}(j)~:~\varphi_{F_{\mathrm{in}}^{(1)}}(F_i^{(1)})(i_1+1+j)=0\}.
\end{equation*}
Note that the designated restriction face contains exactly the vertices that lie in the initial facet, whereas $A_i^{(1)}$ contains the vertices not in $F_{\mathrm{in}}^{(1)}$.
Let $i\geq 1$. First note that $A_i\neq \emptyset$ and $R_i\neq \emptyset$. Since the facets are ordered according to the degree lexicographic ordering with respect to $F_0^{(1)}$, the only facets of $\Delta$ containing $A_i^{(1)}$ are those that are larger than $F_i^{(1)}$ and which hence belong to $\Delta_i$. Therefore, $A_i\in\overset{\circ}{\Delta_i}$. Consider a facet $H$ of $
\partial(A_i^{(1)})\ast\langle R_i^{(1)}\rangle$, i.e., $H=A_i^{(1)}\setminus\{F_i^{1}(j)\}\cup R_i^{(1)}$ for some $i_\ell+1\leq j\leq d$ with $\varphi_{F_{0}^{(1)}}(F_i^{(1)})(j)=1$. Then, $H'=H\cup \{F_0^{(1)}(j)\}\in \mathcal{F}(\Diamond(\Gamma_{i_1}))$ and $H'<F_i$, which implies that $H'\notin \Delta_i$. We conclude $H\in \partial\Delta_i$. The operation $\Delta_i\mapsto \Delta_{i}\setminus F_i$ is hence an elementary shelling.

Now let $\ell\geq 2$. We set $\Delta_0=\Delta\setminus \Diamond(\Gamma_{i_1},\ldots,\Gamma_{i_{\ell-1}})$ and $\Delta_i=\Delta_0\setminus \langle F_0^{(\ell)},\ldots, F_{i-1}^{(\ell)}\rangle$ for $1\leq i\leq r_\ell$. It follows from \Cref{lemma:free facet} that the operation $\Delta_0\mapsto \Delta\setminus F_0^{(\ell)}$ is an elementary shelling. To simplify notation we set $\varphi=\varphi_{F_{\mathrm{in}}^{(\ell)}}$ and $F_i=F_i^{(\ell)}$ for $0\leq i\leq r_\ell$. 
For $1\leq i\leq r_\ell$, we let $t(i)$ be the first position at which $F_i$ differs from the initial facet $F_0$, i.e., $t(i)=\min\{0\leq j\leq d-i_\ell-1~:~\varphi(F_i)(j)=1\}+i_\ell+1$. Let $1\leq i\leq r_\ell$ fixed and set $A'=\{i_j~:~1\leq j\leq\ell-1 \text{ and } i_j<i_\ell\}$. 
We distinguish several cases:\\
{\sf Case 1:} $i_1<i_2$\\
For $1\leq i\leq r_\ell$ we define
$$
A_i=A'\cup \{F_i(j)~:~\varphi(F_i)(j)=1\}
$$
and
$$
R_i=(\{0,\ldots,i_\ell-1\}\setminus A') \cup\{v_{i_\ell}\}\cup\{F_i(j)~:~\varphi(F_i)(j)=0\}.
$$
We will show that $A_i$ and $R_i$ satisfy conditions (1)--(3) from \Cref{definition:shelling}. 

{\sf Condition (1):} 
First note  that $A_i\neq \emptyset$ since there always exists $j$ with $\varphi(F_i)(j)=1$ (as $F_i\neq F_0$). Moreover, we have $R_i\neq \emptyset$ since $v_{i_\ell}\in R_i$. Thus, condition (1) holds.

{\sf Condition (2):} 
We need to show that $A_i\in \overset{\circ}{\Delta_i}$. By the argument at the beginning of this proof, we need to prove that $A_i$ is not contained in any facet $H$ of $\Diamond(\Gamma_{i_1},\ldots,\Gamma_{i_{\ell-1}})\cup\langle F_0,\ldots,F_{i-1}\rangle$. Those come in two different types: 
\begin{itemize}
\item[\textbf{Type 1:}] $H\in \{F_0,\ldots,F_{i-1}\}$.
\item[\textbf{Type 2:}] $H\in \Diamond(\Gamma_{i_j})$ for some $1\leq j\leq \ell-1$ with $i_j<i_\ell$.
\end{itemize}
 Since the facets are ordered lexicographically, $F_i$ is the lexicographically smallest facet of $\Diamond(\Gamma_{i_\ell})$ containing $\{F_i(j)~:~\varphi(F_i)(j)=1\}$ and in particular $A_i$. Therefore, $A_i$ cannot be contained in a facet of Type 1. 
If $H\in \Diamond(\Gamma_{i_j})$ is of Type 2, then $i_j\in A'\subseteq A_i$ but $i_j\notin H$ by \Cref{lemma:diamond complexes}, which implies $A_i\not\subseteq H$.

{\sf Condition (3):} We have to show that $\partial A_i\ast\langle R_i\rangle\subseteq \partial\Delta_{i}$. Let $v\in A_i$ and $B=(A_i\setminus\{v\})\cup R_i\in \partial A_i\ast\langle R_i\rangle$ be a facet. If $v=i_j$ for some $i_j\in A'$, then $B=F_i\setminus \{i_j\}\subseteq F_i\setminus \{i_j\}\cup\{v_{i_j}\}$. \Cref{lemma:diamond complexes} implies that the latter is a facet in $\Diamond(\Gamma_{i_j})$, which by construction does not lie in $\Delta_i$. It follows that $B\subseteq \partial\Delta_i$. Suppose $v=F_i(j)$ for some $i_\ell+1\leq j\leq d$ with $\varphi(F_i)(j-i_\ell-1)=1$. In this case, we have $B=F_i\setminus \{F_i(j)\}\subseteq (F_i\setminus \{F_i(j)\})\cup\{F_0(j)\}$. By \Cref{lemma:diamond complexes}, the latter is a facet in $\Diamond(\Gamma_{i_\ell})$ of a smaller degree than $F_i$ and that is hence lexicographically smaller than $F_i$, meaning that $(F_i\setminus \{F_i(j)\})\cup\{F_0(j)\}\notin \Delta_i$. Hence, $B\in \partial\Delta_i$. 

{\sf Case 2:} $i_1>i_2$.\\
We have different subcases: 
\begin{itemize}
\item[(i)] $t(i)\leq i_1$, or
\item[(ii)] $t(i)>i_1$.
\end{itemize} 
In case (i), we define $A_i$ and $R_i$ as in Case 1 and show that conditions (1)--(3) from \Cref{definition:shelling} hold. The arguments are verbatim the same as in Case 1. Only for condition 2 we need to verify that $A_i$ does not lie in a facet of $\Diamond(\Gamma_{i_1})$. First assume $t(i)< i_1$. Then there exists $0\leq j< i_1-i_\ell-1$ with $\varphi(F_i)(j)=1$. Using the definition of $F_0$, we infer $F_i(j)=v_{j+i_\ell+1}\in A_i$. It follows from  \Cref{lemma:diamond complexes} that $v_j\notin\Diamond(\Gamma_{i_1})$  as $0\leq j\leq i_1-1$. We conclude that $A_i$ cannot be contained in a facet of $\Diamond(\Gamma_{i_1})$. Next assume $t(i)=i_1$. In this case, it follows that $F_i(i_1)=i_1\in A_i$. As $i_1\notin \Diamond(\Gamma_{i_1})$, the claim follows. 

If we are in case (ii) we define
$$
A_i=A'\cup \{v_{i_\ell}\}\cup\{F_i(j)~:~\varphi(F_i)(j)=1\}
$$
and
$$
R_i=(\{0,\ldots,i_\ell-1\}\setminus A') \cup\{F_i(j)~:~\varphi(F_i)(j)=0\}.
$$
Once more, we need to verify conditions (1)--(3).

{\sf Condition 1:}
$t(i)>i_1$ implies that $F_i(i_1)=F_0(i_1)=v_{i_1}$ and hence $\{F_i(j)~:~\varphi(F_i)(j)=0\}\neq \emptyset$. It follows that $R_i\neq \emptyset$. We also have $A_i\neq \emptyset$ since $v_{i_\ell}\in A_i$.

{\sf Condition 2:} The same arguments as above show that $A_i$ is not contained in a facet of type 1 or 2. Since $v_{i_\ell}\in A_i$ but $v_{i_\ell}\notin \Diamond(\Gamma_{i_1})$ (by \Cref{lemma:diamond complexes}), the face $A_i$ cannot be contained in a facet of $\Diamond(\Gamma_{i_1})$.

{\sf Condition 3:} The same arguments as above show that $F_i\setminus \{i_j\}\subseteq \partial\Delta_i$ for $2\leq j\leq \ell-1$ and $F_i\setminus\{F_i(j)\}\subseteq \partial\Delta_i$ for $i_\ell+1\leq j\leq d$ with $\varphi(F_i)(j)=1$. It only remains to show that $F_i\setminus \{v_{i_\ell}\}=(A_i\setminus \{v_{i_\ell}\})\cup R_i\subseteq \partial\Delta_i$. As $t(i)>i_1$, it holds that 
$$
\{0,\ldots,i_\ell-1,i_\ell+1,\ldots,i_1-1,v_{i_1}\}\subseteq F_i\setminus \{v_{i_\ell}\}\subseteq(F_i\setminus\{v_{i_\ell}\})\cup\{i_\ell\}.
$$
\Cref{lemma:diamond complexes} then implies that $(F_i\setminus\{v_{i_\ell}\})\cup\{i_\ell\}\in \Diamond(\Gamma_{i_1})$. The claim follows. 
\end{proof}

%
%
%
 \begin{example}
 \Cref{right_shelling} shows the shelling order from \Cref{theorem:shelling} for $\Diamond(\Gamma_0,\Gamma_1,\Gamma_2)$ within a $2$-ball. First $\Diamond(\Gamma_1)$ is shelled, since it contains the first free facet. Afterwards, one continues with $\Diamond(\Gamma_0)$ and $\Diamond(\Gamma_1)$. 
 	\begin{figure}[h]
 		\centering
 		\includegraphics[scale=0.7]{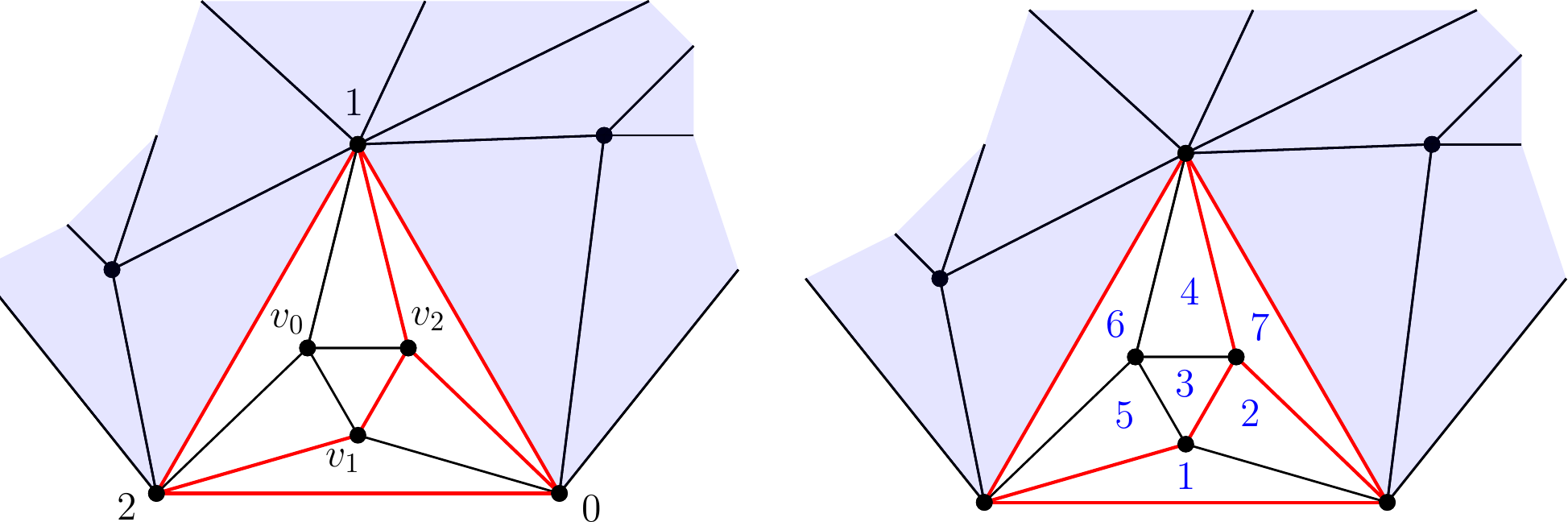}
 		\caption{The blue numbers on the right represent the shelling order on $\Diamond(\Gamma_0,\Gamma_1,\Gamma_2)$. The red lines separate the three complexes $\Diamond(\Gamma_i)$.}
 		\label{right_shelling}
 	\end{figure} 
 \end{example}

\subsection{Proof of Theorem \ref{mainresult}} \label{Sect:ProofMainResult}
\begin{figure}[h]
		\centering
		\includegraphics[scale=0.7]{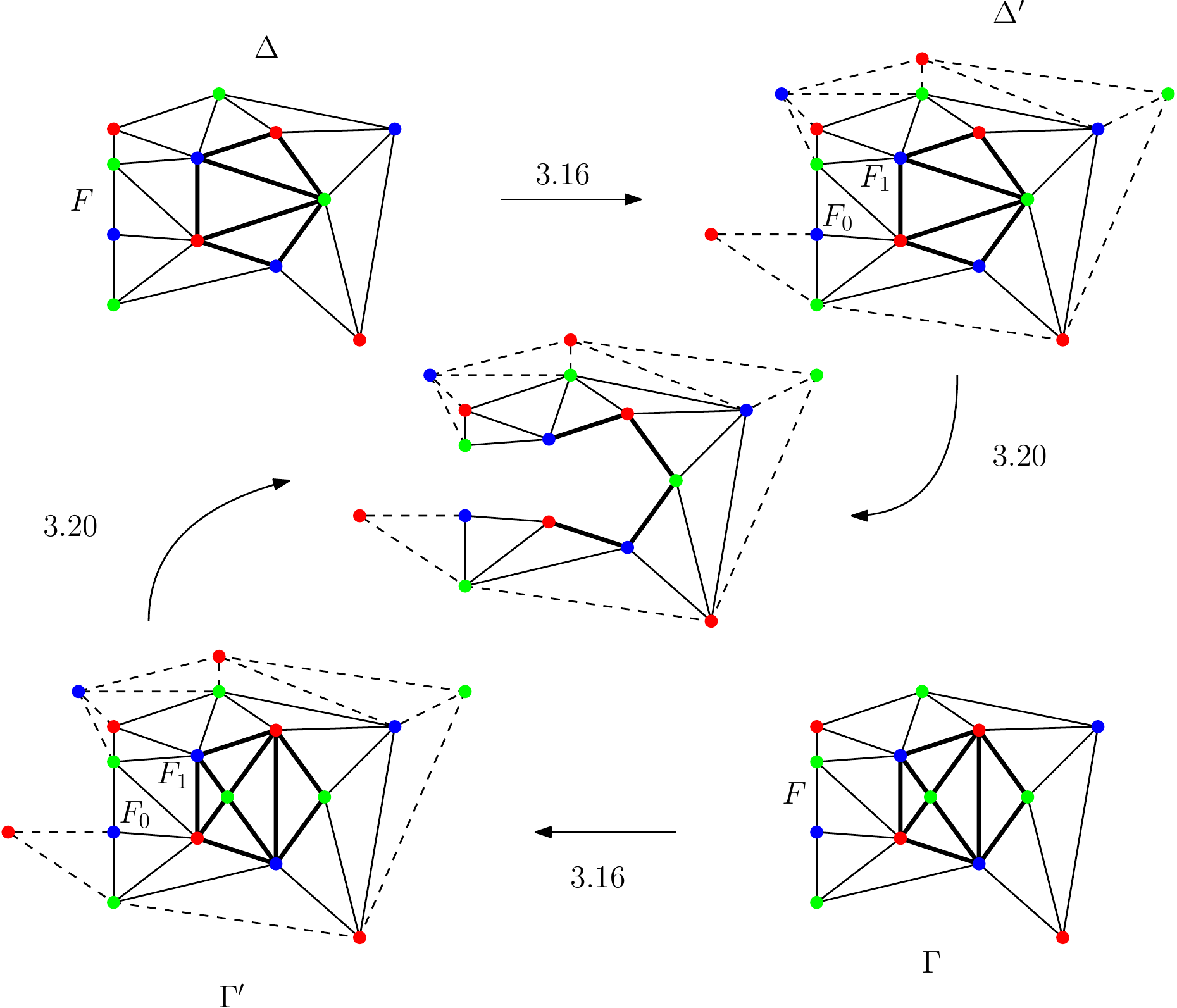}
		\caption{An illustration of Theorem \ref{mainresult} for $\Gamma=\chi^*_D(\Phi)$. The labelling of the faces follows the one in the proof.}
		\label{collar_figure}
	\end{figure}

Combining the results from Sections \ref{sect:same boundary} -- \ref{sect:collar and shelling} we can finally provide the proof of \Cref{mainresult}. 
An illustration of the proof is given in \Cref{collar_figure}.\\

{\sf Proof of \Cref{mainresult}:} 
	Let $\Delta$ and $\Gamma$ be balanced combinatorial $d$-manifolds with boundary. If $\Delta\bsh\Gamma$, then clearly $\Delta$ and $\Gamma$ are PL homeomorphic.\\
	Assume that $\Delta$ and $\Gamma$ are PL homeomorphic. We need to show that $\Delta$ and $\Gamma$ are related by a sequence of shellings and their inverses, preserving balancedness in each step. 
	Using \Cref{bound_same}, we can assume that the boundaries of $\Delta$ and $\Gamma$ are isomorphic and that this isomorphism respects the coloring. It follows from \Cref{cross_int} that there exists a sequence of basic cross-flips, connecting $\Delta$ and $\Gamma$. W.l.o.g. we can assume that $\Delta$ and $\Gamma$ differ by a single cross-flip, i.e., $\Gamma=\chi^*_D(\Delta)$ with $D=\Diamond(\Phi)$, for some subcomplex $\Phi\subseteq\partial\sigma^{d+1}$.
	We can also assume that $\Delta$ (and $\Gamma$) are connected, which implies that they are strongly connected and that there exists a sequence of facets $F_0,\dots,F_m\in\Delta$ such that $F_{i}$ and $F_{i+1}$ intersect in a common face of dimension $d-1$ (for $0\leq i\leq m-1$) and $F_0$ and $F_m$ intersect $\partial \Delta$ and $D$, respectively, in a face of dimension $d-1$ and $d$, respectively. Choosing such a sequence with minimal $m$, we can assure that $F_0,\ldots,F_{m-1},F_m'$~--~ where $F_m'$ is the unique face in $\C_d\setminus D$ that intersects $F_m$ in $F_m\cap \partial D$~--~ is a sequence of facets in $\Gamma$ with the same properties. We now proceed by induction on $m$. Let $F=\partial \Delta\cap F_0$. By assumption we have that $\dim F=d-1$. Applying \Cref{collar} to $\Delta$, we can construct a sequence of inverse shellings that transforms $\Delta$ into a balanced manifold $\Delta'$ such that $\partial\Delta'\cap\partial\Delta=\left\langle F\right\rangle$. Since $\Delta$ and $\Gamma$ have isomorphic boundaries, we can apply the same sequence of inverse shellings to $\Gamma$ in order to obtain a balanced manifold $\Gamma'$, whose boundary is isomorphic to $\partial\Delta'$. 
	 As $\Delta$ and $\Gamma$ are induced subcomplexes of $\Delta'$ and $\Gamma'$, respectively, we can apply the cross-flip $\chi^*_D$ and $\chi^*_{\C_d\setminus D}$ to $\Delta'$ and $\Gamma'$, respectively. In particular, we have $\Delta'\setminus D=\Gamma'\setminus (\C_d\setminus D)$. 
	If $m=0$, then we have $D\cap \partial \Delta' =F$. \Cref{theorem:shelling} further yields a sequence of shellings transforming $\Delta'$ into $\Delta'\setminus D$ and a sequence of inverse shellings from $\Delta'\setminus D=\Gamma'\setminus (\C_{d}\setminus D)$ to $\Gamma'$:
	\begin{equation}\label{eq:Case2}
	\Delta'\mapsto \Delta_1\mapsto \ldots \mapsto \Delta_s=\Delta'\setminus D=\Gamma'\setminus(\C_d\setminus D)\mapsto \Gamma_r\mapsto \ldots\mapsto\Gamma_1\mapsto \Gamma'.
	\end{equation}
	Since $\Delta'$ and $\Gamma'$ are both balanced and every $\Delta_i$ and $\Gamma_i$ is a subcomplex of $\Delta'$ and $\Gamma'$, respectively, every intermediate step in \eqref{eq:Case2} is balanced. We conclude that $\Delta'\bsh\Gamma'$ and hence, by construction of $\Delta'$ and $\Gamma'$, we also have $\Delta\bsh \Gamma$.\\
	Now let $m\geq 1$. We set $\Delta''=\Delta'\setminus F_0$ and $\Gamma''=\Gamma'\setminus F_0$. As $m$ is minimal, we infer that $D\subseteq \Delta''$ and $\C_d\setminus D\subseteq \Gamma''$. Moreover, since $\Delta$ and $\Gamma$ are induced subcomplexes of $\Delta'$ and $\Gamma'$, respectively, it follows that we can apply the cross-flip $\chi^*_D$ and $\chi^*_{\C_d\setminus D}$ to $\Delta''$ and $\Gamma''$, respectively. 
	Moreover, we have that $\partial\Delta''\cong \partial\Gamma''$ and 
	$$
	\Gamma''=\chi^*_D(\Delta'')=\chi^*_D(\Delta')\setminus F_0=\Gamma'\setminus D.
	$$
	Applying the induction hypothesis to $\Delta''$ and $\Gamma''$, yields a sequence of shellings and inverse shellings from $\Delta''$ to $\Gamma''$. It hence follows that $\Delta'\bsh\Gamma'$ and thus also $\Delta\bsh \Gamma$.
	\qed
	
	\section{Combinatorics of basic and reducible cross-flips}\label{sect:BasicCrossflips}
	The aim of this section is threefold. In the first and second part, we will provide a solution to \Cref{problem:cross-flips} by determining the number of combinatorially different basic cross-flips (\Cref{theorem:NumberCrossflips}) and by computing the $h$-vectors of diamond complexes. However, from experiments it is clear that not all of those moves are really needed to transform two closed balanced PL homeomorphic manifolds into each other. A simple reason for this might be that one can write a certain cross-flip as a combination of others. In the second part, we will provide a set of such moves, which will be called \emph{reducible} moves.
	
	\subsection{Counting basic cross-flips}\label{sect: number cross-flips}
\begin{figure}[h]
	\centering
	\includegraphics[scale=0.8]{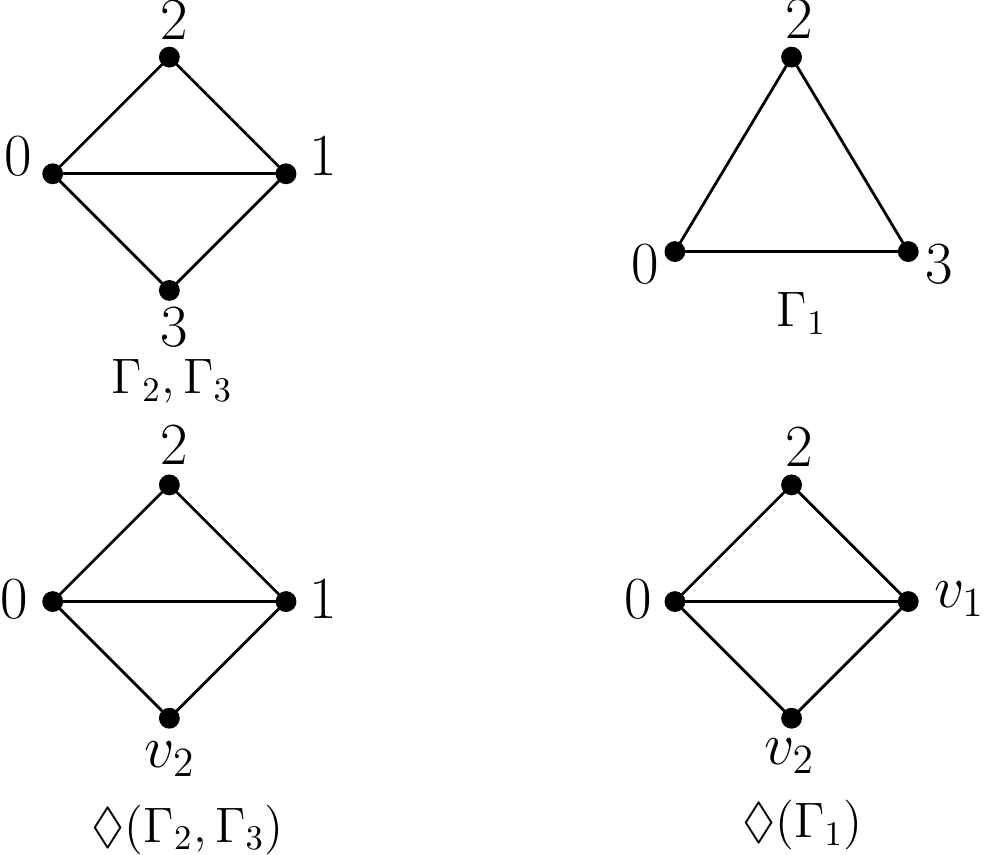}
	\caption{The isomorphism $\Diamond(\Gamma_2,\Gamma_3)\cong\Diamond(\Gamma_1)$ for $d=2$.}
	\label{isomorphic}
\end{figure}

	As a warm-up, we consider the possible basic cross-flips in dimension $2$. \Cref{cross} depicts all combinatorially different (not necessarily basic) cross-flips. The basic cross-flips can be seen in the first two rows of this figure. More precisely, the left picture in the first row shows the cross-flip that exchanges $\Diamond(\Gamma_2)\cong \Diamond(\Gamma_3)$ with $\Diamond(\Gamma_ 0,\Gamma_1,\Gamma_3)\cong \Diamond(\Gamma_0,\Gamma_1,\Gamma_2)$ (and its reverse). It follows from \Cref{lemma:diamond complexes} that $\Diamond(\Gamma_{d})$ and $\Diamond(\Gamma_{d+1})$ are always just $d$-simplices and in particular isomorphic. The isomorphism $\Diamond(\Gamma_ 0,\Gamma_1,\Gamma_3)\cong \Diamond(\Gamma_0,\Gamma_1,\Gamma_2)$ can also be generalized appropriately to higher dimensions and more complicated subcomplexes (see \Cref{isomorphism}). In the right picture on the first line of \Cref{cross} one sees the cross-flip, that removes $\Diamond(\Gamma_1)\cong\Diamond(\Gamma_2,\Gamma_3)$ and adds $\Diamond(\Gamma_0,\Gamma_2,\Gamma_3)\cong\Diamond(\Gamma_0,\Gamma_1)$ (and vice versa). Those isomorphisms will also be instances of \Cref{isomorphism} below. The isomorphism  $\Diamond(\Gamma_1)\cong\Diamond(\Gamma_2,\Gamma_3)$ is also shown in \Cref{isomorphic}. The left picture in the second row of \Cref{cross} shows the cross-flip (and its inverse) substituting  $\Diamond(\Gamma_0)$ with $\Diamond(\Gamma_1,\Gamma_2,\Gamma_3)$, which turns out to be isomorphic to $\Diamond(\Gamma_0)$. We will refer to this move as the \emph{trivial move} and such a move exists in each dimension. Finally, the right picture in the second row depicts the cross-flip interchanging $\Diamond(\Gamma_1,\Gamma_3)\cong \Diamond(\Gamma_1,\Gamma_2)$ and $\Diamond(\Gamma_0,\Gamma_2)\cong\Diamond(\Gamma_0,\Gamma_3)$ and vice versa. Again, the isomorphisms will follow from \Cref{isomorphism}.

We state a first simple observation.
\begin{lemma}\label{isomorphism1}
Let $0\leq k\leq d$. Then 
$$
\Diamond(\Gamma_k)\cong \bigcup_{i=k+1}^{d+1}\Diamond(\Gamma_{i}).
$$
\end{lemma}	
\begin{proof}
The claimed isomorphism is provided by switching vertices $k$ and $v_k$. The statement is then immediate since on the one hand, by \Cref{lemma:diamond complexes}, $F\in\Diamond(\Gamma_k)$ is a facet if and only $\{0,\ldots,k-1,v_k\}\subseteq F$, and on the other hand, $G\in \bigcup_{i=k+1}^{d+1}\Diamond(\Gamma_{i})$ is a facet if and only if $\{0,\ldots,k\}\subseteq G$. 
\end{proof}

The previous lemma in particular provides an explanation for the isomorphism of the trivial move and also for $\Diamond(\Gamma_1)\cong \Diamond(\Gamma_2,\Gamma_3)$, as seen above. However, it is a priori not clear that 
this isomorphism is somehow preserved, when we add further subcomplexes, e.g., we cannot use it to explain the isomorphism $\Diamond(\Gamma_0,\Gamma_1)\cong \Diamond(\Gamma_0,\Gamma_2,\Gamma_3)$. This is the content of the next lemma.

\begin{lemma}\label{isomorphism}
Let $k\in \mathbb{N}$, $0\leq \ell\leq d$ and $0\leq i_1<\cdots <i_k<\ell$. Then,
$$
\Diamond(\Gamma_{i_1},\ldots,\Gamma_{i_k},\Gamma_\ell)\cong\Diamond(\Gamma_{i_1},\ldots,\Gamma_{i_k},\Gamma_{\ell+1},\ldots,\Gamma_{d+1}).
$$
\end{lemma}

\begin{proof}
For $k=0$, the statement is \Cref{isomorphism1}. Now assume $k\geq 1$. We set $D=\Diamond(\Gamma_{i_1},\ldots,\Gamma_{i_k},\Gamma_\ell)$ and $D'= \Diamond(\Gamma_{i_1},\ldots,\Gamma_{i_k},\Gamma_{\ell+1},\ldots,\Gamma_{d+1})$. Let further $\tilde{D}=\Diamond(\Gamma_{i_1},\ldots,\Gamma_{i_k})$. As $k\geq 1$, we deduce from \Cref{lemma:diamond complexes} that $V(D)=V(D')$ and $\{\ell,v_\ell\}\subseteq V(D)$. Let $\psi_\ell:\;V(D)\to V(D')$ be the map that switches vertices $\ell$ and $v_\ell$, i.e.,
$$
\psi_\ell(v)=\begin{cases}
v \quad&\mbox{ if }v\notin \{\ell,v_\ell\},\\
\ell \quad &\mbox{ if } v=v_\ell,\\
v_\ell \quad &\mbox{ if } v=\ell.
\end{cases}
$$
We claim that $\psi_\ell$ induces a simplicial isomorphism between $D$ and $D'$. As $\ell>i_k$ it follows from \Cref{lemma:diamond complexes} that $G\cup \{\ell\}$ is a facet of $\tilde{D}$ if and only $G\cup\{v_\ell\}$ is a facet of $\tilde{D}$. Hence, $F\in \mathcal{F}(D)\cap \mathcal{F}(\tilde{D})$ if and only $\psi_\ell(F)\in \mathcal{F}(D')\cap \mathcal{F}(\tilde{D})$. Moreover, $F$ is a facet in $\Diamond(\Gamma_{\ell})$ if and only if $\{0,\ldots,\ell-1,v_\ell\}\subseteq F$, if and only if $\{0,\ldots,\ell-1,\ell\}\subseteq \psi_\ell(F)$, which is the case if and only if $\psi_\ell(F)$ is a facet of $\Diamond(\Gamma_{\ell+1},\ldots,\Gamma_d)$ by \Cref{lemma:diamond complexes}. The claim follows.
\end{proof}

We state the main result of this section which provides a solution to \Cref{problem:cross-flips}.
\begin{theorem}\label{theorem:NumberCrossflips}
There are $2^{d+1}-1$ combinatorially different basic cross-flips in dimension $d$ 
\end{theorem}
\begin{proof}
It follows from \Cref{lemma:diamond complexes} that $f_{d}(\Diamond(\Gamma_{\ell}))=2^{d-\ell}$ for $0\leq \ell \leq d$ and $f_d(\Diamond(\Gamma_{d+1}))=1$. In particular, if $0\leq i_1<i_2<\cdots <i_k\leq d$, then $f_d(\Diamond(\Gamma_{i_1},\ldots,\Gamma_{i_k}))=\sum_{\ell=1}^{k}2^{d-i_\ell}$. As the representation of $f_d(\Diamond(\Gamma_{i_1},\ldots,\Gamma_{i_k}))$ as a sum of different powers of $2$ is clearly unique, it follows that 
$$
f_d(\Diamond(\Gamma_{i_1},\ldots,\Gamma_{i_k}))\neq f_d(\Diamond(\Gamma_{j_1},\ldots,\Gamma_{j_s})),
$$
for distinct sequences $0\leq i_1<\cdots <i_k\leq d$ and $0\leq j_1<\cdots <j_s\leq d$. This implies that $\Diamond(\Gamma_{i_1},\ldots,\Gamma_{i_k})$ and $\Diamond(\Gamma_{j_1},\ldots,\Gamma_{j_s})$ cannot be isomorphic and the only isomorphisms we have to consider are those given by \Cref{isomorphism}. Therefore, the number of combinatorially different basic cross-flips is given by the number of non-empty subsets of $\{0,\ldots,d\}$ and the claim follows. 
\end{proof}

\subsection{Face numbers of basic cross-flips}\label{sect:numberCrossflips}
In this section, our aim is to compute the face numbers of the diamond complexes that describe basic cross-flips. This will be done by showing that the degree lexicographic shelling order from \Cref{theorem:shelling} for $(\Delta,\Delta\setminus\Diamond(\Gamma))$ also provides a shelling order for $\Diamond(\Gamma)$. 

\begin{proposition}\label{prop:shell}
 Let $\Gamma=\langle \Gamma_{i_1},\ldots,\Gamma_{i_k}\rangle$, where $0\leq i_1<i_2<i_3<\cdots<i_k\leq d+1$. For $1\leq \ell \leq k$, let $r_\ell=f_d(\Diamond(\Gamma_{i_\ell}))-1$ and let 
$$
F_0^{(\ell)}=F_{\mathrm{in}}^{(\ell)},F_1^{(\ell)},\ldots,F_{r_\ell}^{(\ell)}
$$ 
be the ordering of the facets of $\Diamond(\Gamma_{i_\ell})$ according to the degree lexicographic ordering with respect to $F_{\mathrm{in}}^{(\ell)}$. Then
$$
F_{0}^{(1)},F_1^{(1)},\ldots,F_{r_1}^{(1)}, \ldots, F_{0}^{(k)},F_1^{(k)},\ldots,F_{r_k}^{(k)}
$$
is a shelling order for $\Diamond(\Gamma)$. 
\end{proposition}

\begin{proof}
Let $F_i^{(\ell)}$ be a facet. We claim that
$$
R_i^{(\ell)}=\{i_1,\ldots,i_{\ell-1}\}\cup\{F_i^{(\ell)}(j)~:~0\leq j\leq d-i_\ell-1\mbox{ such that } \varphi_{F_0^{(\ell)}}(F_i^{(\ell)})(j)=1\}
$$
is the restriction face of $F_i^{(\ell)}$, i.e., the unique minimal face of $F_i^{(\ell)}$ not lying in $\Delta_i^{(\ell)}=\bigcup_{j=1}^{\ell-1}\Diamond(\Gamma_{i_j})\cup\langle F_0^{(\ell)},\ldots,F_{i-1}^{(\ell)}\rangle$. 
Let $H\in \langle F_i^{(\ell)}\rangle\setminus \Delta_i^{(\ell)}$. As $H\not\subseteq F_j^{(\ell)}$ for $0\leq j\leq i-1$, we must have $\{F_i^{(\ell)}(j)~:~\varphi_{F_0^{(\ell)}}(F_i^{(\ell)})(j)=1\}\subseteq H$. Otherwise, $H$ would be contained in a degree lexicographically smaller facet than $F_i^{(\ell)}$. Moreover, $H\not\subseteq \Diamond(\Gamma_{i_j})$ for $1\leq j\leq \ell-1$ implies that $i_j\in H$ for $1\leq H\leq \ell-1$. We conclude that $R_i^{(\ell)}\subseteq H$. The same arguments also show that indeed $R_i^{(\ell)}\notin \Delta_i^{(\ell)}$. The claim follows.
\end{proof}

We use the shelling of \Cref{prop:shell} to compute the $h$-vectors of diamond complexes. In the following, we set $\binom{n}{k}=0$ if $n<k$. 

%
	\begin{proposition}\label{prop:hVectorDiamond}
	Let $0\leq i_1<\cdots <i_k\leq d+1$ and let $\Diamond(\Gamma)=\Diamond(\Gamma_{i_1},\ldots,\Gamma_{i_k})$. Then, for $0\leq \ell\leq d$
		$$h_\ell(\Diamond(\Gamma))=\binom{d-i_1}{ \ell}+\binom{d-i_2}{\ell-1}+\cdots+\binom{d-i_k}{\ell-k+1}.$$
	\end{proposition}
	
	The previous proposition in particular implies that $h_\ell(\Diamond(\Gamma))=0$ for $\ell>d-i_1$. Also note that for $k=1$, the statement follows directly from \Cref{lemma:diamond complexes} combined with the fact that $h_\ell(\C_d)=\binom{d+1}{\ell}$.
	
	\begin{proof}
	By \eqref{eq:shelling},  to compute $h_\ell(\Diamond(\Gamma))$ we need to count restriction faces of size $\ell$. It follows from the proof of \Cref{prop:shell} that for $1\leq j\leq k$ the facets in $\Diamond(\Gamma_{i_j})$, whose restriction faces  are of size $\ell$ are those that differ in $\ell-(j-1)$ positions from the initial facet $F_0^{(j)}$. Since there are $\binom{d-i_j}{\ell-j+1}$ such facets, one for each $(\ell-j+1)$-subset of $\{0,\ldots,d-i_j\}$, the claim follows.
	\end{proof}
Using the following lemma, we can control the face numbers when applying a cross-flip.
	
		\begin{lemma}\label{lemma:h}
			Let $D$ be a shellable and co-shellable subcomplex of $\mathcal{C}_d$ that is a $d$-ball. Then, for $0\leq i\leq d+1$
			$$h_i(D)+h_{d+1-i}(\mathcal{C}_d\setminus D)=\binom{d+1}{i}.$$
		\end{lemma}
		\begin{proof}
			It is known that 
			$$h_i(\mathcal{C}_d)=\binom{d+1}{i}$$
			for all $0\leq i\leq d$. Since $D$ is both shellable and co-shellable, any shelling order on $D$ can be extended to one of $\mathcal{C}_d$ by adding a reverse order on $\mathcal{C}_d\setminus D$. 
			Reversing the order on $\C_d\setminus D$ has the effect that a facet $F$, whose restriction face $R$ was of size $i$ before, now has the restriction face $F\setminus R$, which is of size $d+1-i$.
			The claim now follows from \eqref{eq:shelling}.
		\end{proof}
		
		\begin{example}
	We use \Cref{lemma:h} to compute the $h$-vector of a cross-polytopal stacked $d$-sphere. A \emph{cross-polytopal stacked $(d-1)$-sphere} $\mathcal{ST}^\times\left( n,d\right)$ on $n$ vertices is defined to be the connected sum of $\left( \frac{n}{d+1}-1\right)$ copies of $\mathcal{C}_d$. We can build $\mathcal{ST}^\times\left( n,d\right)$ from $\mathcal{C}_d$, by replacing  a facet $F$ with its complement in $\mathcal{C}_d$  and repeating this process $\left(\frac{n}{d+1}-2\right)$-times. In other words, we apply the basic cross-flip $\chi^*_{\Diamond(\Gamma_d)}$ $\left(\frac{n}{d+1}-2\right)$-times to $\C_d$.  By \Cref{lemma:h}
	$$h_{d+1-i}(\mathcal{C}_d\setminus F)=\binom{d+1}{i}$$
	for $0\leq i\leq d$ and $h_{d+1}(\C_d\setminus F)=0$. 
	Using a reversed shelling of $\C_d\setminus F$, we can successively transform $\C_d$ into $\mathcal{ST}^\times\left( n,d\right)$ and obtain
	$$h_i(\mathcal{ST}^\times\left( n,d\right) )=\binom{d+1}{i}+\left( \frac{n}{d+1}-2\right) \binom{d+1}{i}=\left( \frac{n}{d+1}-1\right) \binom{d+1}{i}$$
	for $0<i<d$ and $h_0(\mathcal{ST}^\times\left( n,d\right))=h_{d+1}(\mathcal{ST}^\times\left( n,d\right))=1$. We want to remark that the $h$-vector of $\mathcal{ST}^\times\left( n,d\right)$ was well-known before and can be computed directly using that the $h$-vector is additive with respect to connected sum. However, we included this example as a nice application of \Cref{lemma:h}.
\end{example}

\subsection{Reducible cross-flips}
In \Cref{sect: number cross-flips} we have seen that there are $2^{d+1}-1$ combinatorially different cross-flips in dimension $d$. However, if we want to relate two balanced closed PL homeomorphic manifolds using basic cross-flips as in \Cref{ikn cross}, then it is conceivable that not all of these are really needed. Indeed, it was shown in \cite{Murai:Suzuki} that for $2$-dimensional spheres, none of which is the boundary of the cross-polytope, it is sufficient to use the pentagon move and its inverse, which are moves on the right in the second line of \Cref{cross}. This is the motivation for this section. 

\begin{definition}
Let $\mathrm{BC}_d$ the set of basic cross-flips in dimension $d$. 
We call a set $R\subseteq\mathrm{BC}_d$ of basic cross-flips \emph{reducible} if every cross-flip in $R$ can be expressed as a combination of cross-flips in $\mathrm{BC}_d\setminus R$. 
\end{definition}
Note that if $R$ is a set of reducible cross-flips, then, in particular, the basic cross-flips contained in $\BC_d\setminus R$ suffice to relate any two balanced closed PL homeomorphic manifolds.
By \Cref{isomorphism} we also know that for any basic cross-flip, there is a representative of the form $\chi^*_{\Diamond(\Gamma)}$, where $\Gamma=\langle \Gamma_{i_1},\ldots,\Gamma_{i_k}\rangle$, and $0\leq i_1<\cdots <i_k\leq d$. In the following, we will always use this representative. 
We need to introduce some further notation. For a set $I\subseteq\{0,\ldots,d\}$ and $a\in \ZZ$, we write $I+a$ for the set $\{i+a~:~i\in I\}$. Moreover, given $I=\{i_1,\ldots,i_k\}$ we use $\Gamma_I$ to denote the $d$-ball $\langle\Gamma_{i_1},\ldots,\Gamma_{i_k}\rangle$. 
Since the dimension of $\Gamma_I$ is not clear from the notation a priori, we will add a superscript and write $\Gamma_I^d$ from now on. Similarly, we write $\Diamond^d$ for the diamond operation in dimension $d$. Those distinctions will be important in the rest of this section.

\begin{figure}[h]
	\centering
	\includegraphics[scale=1.2]{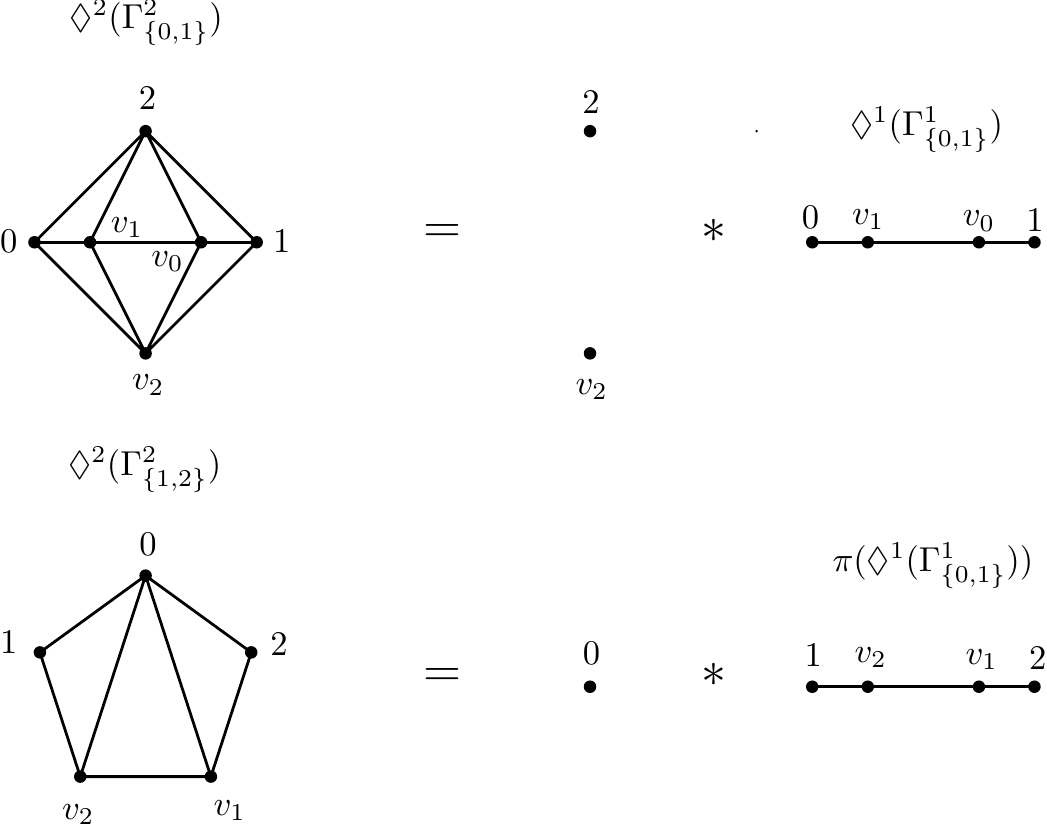}
	\caption{The isomorphism from \Cref{lemma:Suspension and Cone}.}
	\label{cone_susp}
\end{figure}

We have the following observations. An instance is depicted in \Cref{cone_susp}.
\begin{lemma}\label{lemma:Suspension and Cone}
Let $I\subseteq \{0,\ldots,d\}$.
\begin{itemize}
\item[(1)] If $d\notin I$, then $\Diamond^d(\Gamma^d_I)=\{d,v_d\}\ast\Diamond^{d-1}(\Gamma_I^{d-1})$.
\item[(2)] If $0\notin I$, then $\Diamond^d(\Gamma^d_I)=\{0\}\ast \pi(\Diamond^{d-1}(\Gamma_{I-1}^{d-1}))$, where 
\begin{align*}
\pi:\{0,\ldots,d-1\}\cup\{v_0,\ldots,v_{d-1}\}&\to \{1,\ldots,d\}\cup\{v_1,\ldots,v_d\}\\
i&\mapsto i+1 \\
v_i&\mapsto v_{i+1}.
\end{align*}
\end{itemize}
\end{lemma}

\begin{proof}
It follows from \Cref{lemma:diamond complexes} that we have $\Diamond^d(\Gamma_j^d)=\{d,v_d\}\ast \Diamond^{d-1}(\Gamma_j^{d-1})$ for $0\leq j< d$. This implies (1).

For (2), first note that $0$ lies in every facet of $\Diamond^d(\Gamma^{d}_j)$ if $j\neq 0$. The statement is then immediate by \Cref{lemma:diamond complexes}.
\end{proof}
Let $\psi:\{d,v_d\}\to \{d,v_d\}$ be the map that interchanges $d$ and $v_d$ and let 
\begin{align*}
\rho:\{0,\ldots,d\}\cup\{v_0,\ldots,v_d\}&\to \{0,\ldots,d\}\cup\{v_0,\ldots,v_d\}\\
 i&\mapsto i-1 \\
v_i &\mapsto v_{i-1},
\end{align*}
where we compute $i-1$ modulo $(d+1)$. 
We also set $\sigma=\psi\circ\rho$. 

\begin{figure}[h]
	\centering
	\includegraphics[scale=1.3]{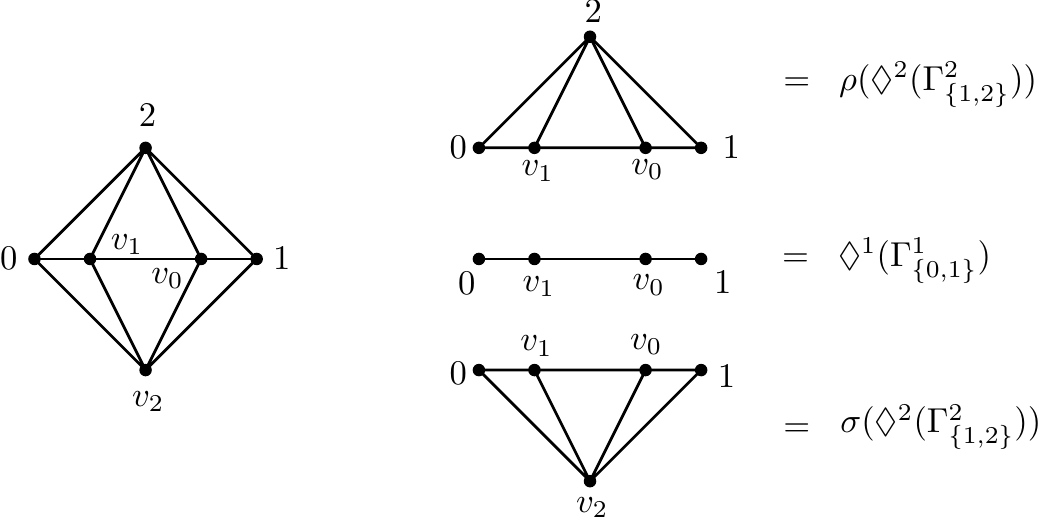}
	\caption{An instance of the decomposition given in \Cref{cor}.}
	\label{sigma_rho}
\end{figure}

As a consequence of part (1) of the previous lemma, we obtain the following decomposition of $\Diamond^d(\Gamma^d_I)$ that will be crucial. Once again we include a 2-dimensional example in \Cref{sigma_rho} for the sake of clarity.
\begin{corollary}\label{cor}
	Let $I\subseteq \{0,\ldots,d\}$ with $d\notin I$. Then,
	\begin{equation}
	 \Diamond^d(\Gamma^d_I)=\rho(\Diamond^d(\Gamma^d_{I+1}))\cup\sigma(\Diamond^d(\Gamma^d_{I+1})),
	\end{equation}
	where 
	$$
	\rho(\Diamond^d(\Gamma^d_{I+1}))\cap\sigma(\Diamond^d(\Gamma^d_{I+1}))=\Diamond^{d-1}(\Gamma^{d-1}_I).
	$$
	Moreover, $\rho(\Diamond^d(\Gamma^d_{I+1}))$ and $\sigma(\Diamond^d(\Gamma^d_{I+1}))$ are induced subcomplexes of $\Diamond^d(\Gamma^d_I)$.	
\end{corollary}

\begin{proof}
    Since by assumption $d\notin I$ and clearly $0\notin I+1$, by \Cref{lemma:Suspension and Cone} (1) we have 
    \begin{align*}
		\Diamond^d(\Gamma^d_I)=&\{d,v_d\}\ast\Diamond^{d-1}(\Gamma_I^{d-1})\\
		=&\left(\{d\}\ast \Diamond^{d-1}(\Gamma_I^{d-1})\right)\cup \left(\{v_d\}\ast \Diamond^{d-1}(\Gamma_I^{d-1})\right).
		\end{align*}
		Both complexes appearing in this decomposition are clearly induced and their intersection is $\Diamond^{d-1}(\Gamma^{d-1}_I)$. 
		Using \Cref{lemma:diamond complexes} one easily verifies that $$
		\{d\}\ast \Diamond^{d-1}(\Gamma_I^{d-1})=\rho(\Diamond^d(\Gamma^d_{I+1}))\quad  \mbox{and}\quad \{v_d\}\ast \Diamond^{d-1}(\Gamma_I^{d-1})=\sigma(\Diamond^d(\Gamma^d_{I+1}).
			$$
\end{proof}

The last ingredient, we need is the following decomposition that follows from part (2) of \Cref{lemma:Suspension and Cone}. 

\begin{corollary}\label{cor2}
Let $I\subseteq \{0,\ldots,d\}$ with $0\in I$. Then,
	\begin{equation}
	 \Diamond^d(\Gamma^d_I)=\Diamond^d(\Gamma^d_{I\setminus \{0\}})\cup\Diamond^d(\Gamma^d_0) 
	\end{equation}
	and 
	$$
	\Diamond^d(\Gamma^d_{I\setminus \{0\}})\cap\Diamond^d(\Gamma^d_0)=\pi(\Diamond^{d-1}(\Gamma^{d-1}_{(I\setminus \{0\})-1})).
	$$
	Moreover, $\Diamond^d(\Gamma^d_0)$ is an induced subcomplex of $\Diamond^d(\Gamma^d_I)$.	
\end{corollary}

\begin{proof}
First note that 
$$
\Diamond^d(\Gamma^d_I)=\Diamond^d(\Gamma^d_{I\setminus \{0\}})\cup\Diamond^d(\Gamma^d_0)
$$
and that $\Diamond^d(\Gamma^d_0)$ is induced. 
Since $0\notin I\setminus \{0\}$, it follows from \Cref{lemma:Suspension and Cone} (2) that 
$$
\Diamond^d(\Gamma^d_{I\setminus \{0\}})=\{0\}\ast \pi(\Diamond^{d-1}(\Gamma_{I\setminus\{0\}-1}^{d-1})).
$$
Combining this with \Cref{lemma:diamond complexes} we obtain
$$
\Diamond^d(\Gamma^d_{I\setminus \{0\}})\cap\Diamond^d(\Gamma^d_0)=\pi(\Diamond^{d-1}(\Gamma^{d-1}_{(I\setminus \{0\})-1})).
$$
\end{proof}

 Our aim is to show, that the set of cross-flips that remove diamond complexes $\Diamond^d(\Gamma^d_I)$ with $d\notin I$ is a set of reducible cross-flips. For this, using \Cref{cor}, we decompose $\Diamond^d(\Gamma^d_I)$ into two copies of $\Diamond^d(\Gamma^d_{I+1})$ that are glued together along a subcomplex of their boundaries, that is itself isomorphic to $\Diamond^{d-1}(\Gamma_I^{d-1})$.
The idea is to first flip one of the subcomplexes that are isomorphic to $\Diamond^d(\Gamma^d_{I+1})$, then to decompose again, using \Cref{cor2}, and then to perform another flip. The next example makes this idea more precise.
 \begin{example}
 Let $d=2$ and let us consider $\Diamond^2(\Gamma^2_1)$. We first decompose $\Diamond^2(\Gamma^2_1)$ into two copies of $\Diamond^2(\Gamma^2_2)$. In this example, this is just two triangles intersecting in an edge (see the left picture in the first row of \Cref{reduction_figure}). We now flip the ``upper'' subcomplex (see the right picture in the first row of \Cref{reduction_figure}). The flipped part is now decomposed again into a part that is isomorphic to $\Diamond^2(\Gamma_{\{1,2\}}^2)$ and a part isomorphic to $\Diamond^2(\Gamma_0^2)$ (shown in white and green respectively in the middle picture of \Cref{reduction_figure}). The second component is grouped together with the untouched copy of $\Diamond^2(\Gamma^2_2)$. The union of those two subcomplexes is isomorphic to $\Diamond^2(\Gamma^2_{\{0,2\}})$ (see left picture at the bottom of \Cref{reduction_figure}) and is substituted by its complement in $\C_2$. We obtain $2$ copies of $\Diamond^2(\Gamma^2_{\{1,2\}})$ glued together along a subcomplex of their boundaries, which is $\C_2\setminus(\Diamond^2(\Gamma^2_1))=\Diamond^2(\Gamma^2_{\{0,1\}})$, where the last equality follows from \Cref{isomorphism}.
\begin{figure}[h]
	\centering
	\includegraphics[scale=0.7]{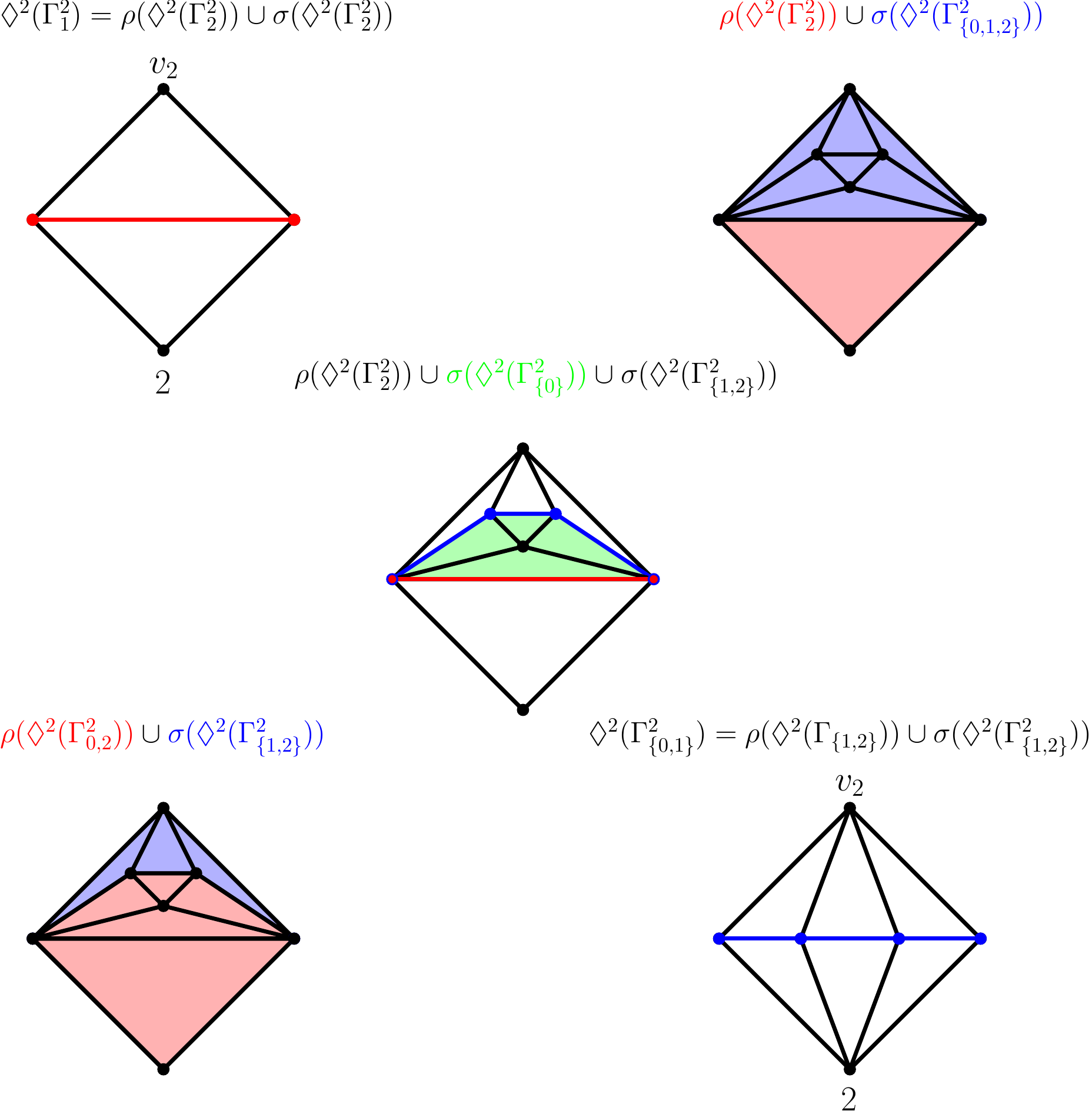}
	\caption{The reducibility of $\chi^*_{\Diamond^2(\Gamma_1^2)}$ and $\chi^*_{\Diamond^2(\Gamma_{\{0,1\}}^2)}$. This particular reduction was already pointed out in \cite[Remark 4.8]{Murai:Suzuki}. }
	\label{reduction_figure}
\end{figure} 
 \end{example}

\begin{remark}\label{remark:complement}
Let $I=\{i_1,\ldots,i_k\}\subseteq \{0,\ldots,d\}$ with $i_1<\cdots<i_k$ and $d\notin I$. The inverse of the cross-flip $\chi^*_{\Diamond^d(\Gamma^d_I)}$ is then given by $\chi^*_{\Diamond^d(\Gamma^d_{I^c})}$, where here and thereafter we set $I^c=\{0,\ldots,d+1\}\setminus I$. However, it follows from \Cref{isomorphism} that the cross-flip $\chi^*_{\Diamond^d(\Gamma^d_J)}$ with 
$$
J=\{0,1,\dots, i_k\}\setminus \{i_1,\dots,i_{k-1}\}
$$
has the same effect. In particular, $d\notin J$. This observation will be useful in the proof of the main theorem of this section, which we now state.
\end{remark}

\begin{theorem}\label{thm:reducible}
The set of basic cross-flips 
$$
\{\chi^*_{\Diamond^d(\Gamma_I^d)}~:~I\subseteq \{0,\ldots, d-1\}\}
$$
is reducible. In other words, the set of basic cross-flips
$$
\{\chi^*_{\Diamond^d(\Gamma_I^d)}~:~I\subseteq \{0,\ldots, d\},\; d\in I\}
$$
is sufficient to relate any two balanced PL homeomorphic closed manifolds.
\end{theorem}

\begin{proof}
Let $I\subseteq \{0,\ldots,d-1\}$. We will show that the basic cross-flip $\chi^*_{\Diamond^d(\Gamma^d_I)}$ can be expressed as the combination of the two basic cross-flips $\chi^*_{\Diamond^d(\Gamma^d_{I+1})}$ and $\chi^*_{\Diamond^d(\Gamma^d_{(I+1)\cup\{0\}})}$. The required statement then follows by applying this substitution process iteratively until $d\in I+1$. \\
By \Cref{cor}, we know that
$$
 \Diamond^d(\Gamma^d_I)=\rho(\Diamond^d(\Gamma^d_{I+1}))\cup\sigma(\Diamond^d(\Gamma^d_{I+1}))
$$
and both subcomplexes on the right are induced. We can hence apply the cross-flip  $\chi^*_{\Diamond^d(\Gamma^d_{I+1})}$, in order to replace one of the components, e.g., the first one, by its complement in $\C_d$. Using that $\rho$ is an isomorphism, that $0\in (I+1)^c$ and \Cref{cor2}, we can compute the complement of $\rho(\Diamond^d(\Gamma^d_{I+1}))$ via
\begin{align*}
\rho(\C_d)\setminus(\rho(\Diamond^d(\Gamma^d_{I+1})))=&\rho(\C_d\setminus\Diamond^d(\Gamma^d_{I+1}))=\rho(\Diamond^d(\Gamma^d_{(I+1)^c}))\\
=&\rho(\Diamond^d(\Gamma^d_{(I+1)^c\setminus \{0\}}))\cup \rho(\Diamond^d(\Gamma^d_{0}))=\rho(\Diamond^d(\Gamma^d_{I^c+1}))\cup \rho(\Diamond^d(\Gamma^d_{0})).
\end{align*}
Note that both complexes, $\sigma(\Diamond^d(\Gamma^d_{I+1}))$ and $\rho(\Diamond^d(\Gamma^d_{0}))$ contain a vertex labeled $v_d$. In $\sigma(\Diamond^d(\Gamma^d_{I+1}))$, this is the vertex corresponding to $0$ in $\Diamond^d(\Gamma^d_{I+1})$, in $\rho(\Diamond^d(\Gamma^d_{0}))$ this is the vertex corresponding to $v_0$ in $\Diamond^d(\Gamma^d_{0})$. We relabel the vertex $v_d\in \rho(\Diamond^d(\Gamma^d_{0}))$ with $\tilde{v}_d$. 
As $\rho(\Diamond^d(\Gamma^d_{0}))=\{\tilde{v}_d\}\ast \C_{d-1}$ and $\sigma(\Diamond^d(\Gamma^d_{I+1}))=\{v_d\}\ast \Diamond^{d-1}(\Gamma^{d-1}_I)$,  their intersection is simply
\begin{align*}
\sigma(\Diamond^d(\Gamma^d_{I+1}))\cap\rho(\Diamond^d(\Gamma^d_{0}))&=\Diamond^{d-1}(\Gamma^{d-1}_I).
\end{align*}
As $0\notin I+1$, it hence follows from \Cref{cor} that 
$$
\sigma(\Diamond^d(\Gamma^d_{I+1}))\cup\rho(\Diamond^d(\Gamma^d_{0}))=\tilde{\sigma}(\Diamond^d(\Gamma_{(I+1)\cup\{0\}}^d)),
$$
where $\tilde{\sigma}$ is the composition of $\sigma$, followed by a relabeling of the vertex $d$ with $\tilde{v}_d$. As by \Cref{cor} and \Cref{cor2} $\sigma(\Diamond^d(\Gamma^d_{I+1}))$ and $\rho(\Diamond^d(\Gamma^d_{0}))$ are induced, so is their union and we can apply the cross-flip $\chi^*_{\Diamond^d(\Gamma^d_{(I+1)\cup\{0\}})}$ substituting $\tilde{\sigma}(\Diamond^d(\Gamma_{(I+1)\cup\{0\}}^d))$ with its complement in $\tilde{\sigma}(\C_d)$, which is given by
$$
\tilde{\sigma}(\C_d)\setminus \tilde{\sigma}(\Diamond^d(\Gamma_{(I+1)\cup\{0\}}^d))=\tilde{\sigma}(\Diamond^d(\Gamma_{((I+1)\cup\{0\})^c}))=\tilde{\sigma}(\Diamond^d(\Gamma_{(I^c+1)\setminus\{0\}})).
$$
It remains to show that the union of the two complexes $\rho(\Diamond^d(\Gamma^d_{I^c+1}))$ and $\tilde{\sigma}(\Diamond^d(\Gamma_{I^c+1}))$ is isomorphic to $\Diamond^d(\Gamma^d_{I^c})$. To see this, first note that by \Cref{remark:complement}, the basic cross-flip $\chi^*_{\Diamond^d(\Gamma^d_{I^c+1})}$ has a representative $\chi^*_{\Diamond(\Gamma_J)}$ with $d\notin J$. The claim then follows from \Cref{cor}. 
\end{proof}

\begin{remark}
As an immediate consequence of \Cref{thm:reducible} it follows that $2^{d}$ basic cross-flips are sufficient to relate any two balanced PL homeomorphic manifolds without boundary, e.g., in dimension $2$ the four cross-flips $\chi^*_{\Diamond^2(\Gamma^2_{2})}$, $\chi^*_{\Diamond^2(\Gamma^2_{\{0,2\}})}$, $\chi^*_{\Diamond^2(\Gamma^2_{\{1,2\}})}$ and $\chi^*_{\Diamond^2(\Gamma^2_{\{0,1,2\}})}$ suffice. In \Cref{cross} the left picture in the first row depicts the interchange of $\Diamond^2(\Gamma^2_{2})$ and $\Diamond^2(\Gamma^2_{\{0,1,2\}})$. The right picture in the second row depicts the interchange of $\Diamond^2(\Gamma^2_{\{1,2\}})$ and $\Diamond^2(\Gamma^2_{\{0,2\}})$.
\end{remark}

\section{Open problems}\label{sect:openProblems}
\subsection{Connecting two balanced manifolds using few flips }
The proof of \Cref{strong thm:crossflips} in \cite{IKN} does not provide information on the number of moves needed to connect two closed balanced combinatorial manifolds that are PL homeomorphic. On the other hand, especially from a computational point of view, an upper bound for the number of operations needed would be of interest. As an example let us assume that we start from a manifold $\Delta$ on $n$ vertices and we perform \emph{every} applicable cross-flip on $\Delta$, obtaining a set $T_1$ of ``target'' manifolds. In the next step we repeat the procedure for each of the objects in $T_1$, and denote with $T_2$ the new targets. We proceed iteratively for a certain number of steps. What is the minimum number of steps that we need to obtain \emph{all} balanced combinatorial manifolds that are PL-homeomorphic to $\Delta$ with $n$ vertices? In a more general way, we formulate the following question.
\begin{question}
	Fix two PL homeomorphic balanced combinatorial manifolds $\Delta$ and $\Gamma$. What is an upper bound (depending on $f_0(\Delta)$ and $f_0(\Gamma)$) for the minimal number of cross-flips needed to connect $\Delta$ and $\Gamma$?
\end{question}
\subsection{The minimal set of sufficient basic cross-flips}
In Section \ref{sect: number cross-flips} we showed that for a fixed dimension $d$, there are precisely $2^{d+1}-1$ combinatorially distinct basic cross-flips, out of which $2^d$ suffice to relate any two balanced PL homeomorphic closed manifolds of dimension $d$. As remarked earlier for the case $d=2$ it is proved in \cite{Murai:Suzuki} that 3 flips actually suffice. Even more surprisingly, there it is shown that it is possible to connect all \emph{but a finite number} of triangulations of  a fixed surface using only 2 flips (see \cite[Theorem 4.3]{Murai:Suzuki}). It is interesting to note that the set of 3 sufficient flips provided by Murai and Suzuki is not contained in the set of 4 flips given in \Cref{thm:reducible}. We hope that our detailed description of these moves will yield an answer to the following question:
\begin{question}
	What is the cardinality of a minimal set of basic cross-flips that suffice to relate any two balanced combinatorial $d$-manifolds that are PL homeomorphic for $d>2$?
\end{question}
\begin{example}
 \Cref{pentagon_red} shows that we can write $\chi^*_{\Diamond^2(\Gamma^2_{\{1,2\}})}$ (or $\chi^*_{\Diamond^2(\Gamma^2_{\{0,2\}})}$ ) as a composition of two flips that belong to the set constructed in \Cref{thm:reducible}. Hence, for $d=2$, the set in \Cref{thm:reducible} can be further reduced to a set of cardinality 3.
\begin{figure}[h]
	\centering
	\includegraphics[scale=0.8]{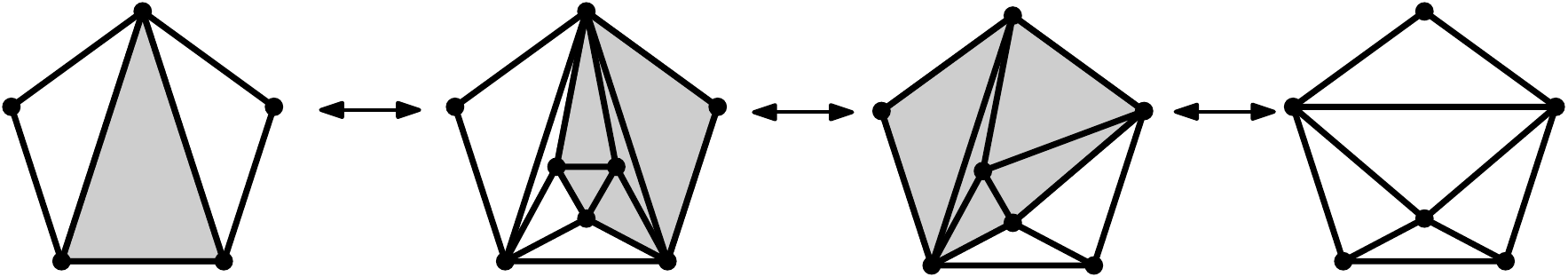}
	\caption{The flip $\chi^*_{\Diamond^2(\Gamma^2_{\{1,2\}})}$ as a composition of $\chi^*_{\Diamond^2(\Gamma^2_{2})}$ and $\chi^*_{\Diamond^2(\Gamma^2_{\{0,2\}})}$. (Reading the figure from right to left clearly shows that $\chi^*_{\Diamond^2(\Gamma^2_{\{0,2\}})}$ can be obtained from $\chi^*_{\Diamond^2(\Gamma^2_{\{0,1,2\}})}$ and $\chi^*_{\Diamond^2(\Gamma^2_{\{1,2\}})}$). }
	\label{pentagon_red}
\end{figure} 
	Moreover, it is not hard to notice by direct inspection that for $d=2$ there are 8 minimal sets of sufficient basic cross-flips, and they all have cardinality 3. More precisely, setting
	\begin{align*}\mathcal{B}=\{&
	\{\{1\},\{0,1\},\{0,2\}\}\\
	&\{\{1\},\{0,1\},\{1,2\}\}\\
	&\{\{1\},\{0,1,2\},\{0,2\}\}\\
	&\{\{1\},\{0,1,2\},\{1,2\}\}\\
	&\{\{2\},\{0,1\},\{0,2\}\}\\
	&\{\{2\},\{0,1\},\{1,2\}\}\\
	&\{\{2\},\{0,1,2\},\{0,2\}\}\\
	&\{\{2\},\{0,1,2\},\{1,2\}\}\},
	\end{align*}
for every $B\in\mathcal{B}$ the set $\{\chi^*_{\Diamond^2(\Gamma^2_{I})}: I\in B\}$ is a set of minimal sufficient cross-flips and there are no other minimal ``generating sets''.
\end{example}
 Since the whole process of reduction seems to encode a notion of dependence between flips, we underline the following property of the set of minimal sufficient cross-flips for $d=2$.
\begin{lemma}
	The set $\mathcal{B}$ is the set of bases of a matroid.
\end{lemma} 
\begin{proof}
	It is straightforward to check that $\mathcal{B}$ satisfies the basis exchange axiom. In particular, the matroid $M(\mathcal{B})$ is isomorphic to a direct sum of 3 uniform matroids on 2 elements of rank 1: $M(\mathcal{B})\cong U_{1,2}\oplus U_{1,2}\oplus U_{1,2}$.
\end{proof}	
The previous lemma, though only providing weak evidence, clearly motivates the following question, which was also raised by Eran Nevo in personal communication:

\begin{question}
	Do all the minimal sets of basic cross-flips that suffice to relate any two PL homeomorphic balanced combinatorial $d$-manifolds form the set of bases of a matroid, for $d\geq 2$? 
\end{question}
Since the previous question appears to be rather ambitious, we also propose the following weaker question. 		
\begin{question}
	Do all the minimal sets of basic cross-flips that suffice to relate any two PL homeomorphic balanced combinatorial $d$-manifolds have the same cardinality?
\end{question}

\section*{Acknowledgements}
We would like to thank Steven Klee and Isabella Novik for their comments on an earlier version of this article that helped to improve it.

\bibliographystyle{alpha}
\bibliography{bibliography}

\begin{thebibliography}{MJKS17}

\bibitem[Ale30]{Ale}
J.~W. Alexander.
\newblock The combinatorial theory of complexes.
\newblock {\em Ann. of Math. (2)}, 31(2):292--320, 1930.

\bibitem[BH93]{BH-book}
W.~Bruns and J.~Herzog.
\newblock {\em Cohen-{M}acaulay rings}, volume~39 of {\em Cambridge Studies in
  Advanced Mathematics}.
\newblock Cambridge University Press, Cambridge, 1993.

\bibitem[BHL03]{Bjoerner:Lutz}
A.~Bj\"{o}rner and F.~H.~Lutz.
\newblock A 16-vertex triangulation of the {P}oincar\'{e} homology 3-sphere and
  non-{PL} spheres with few vertices.
\newblock 01 2003.

\bibitem[Bin83]{Bing}
R.~H. Bing.
\newblock {\em The geometric topology of 3-manifolds}, volume~40 of {\em
  American Mathematical Society Colloquium Publications}.
\newblock American Mathematical Society, Providence, RI, 1983.

\bibitem[Bj{\"o}95]{Bjoerner}
A.~Bj{\"o}rner.
\newblock Topological methods.
\newblock In {\em Handbook of combinatorics, {V}ol.\ 1,\ 2}, pages 1819--1872.
  Elsevier Sci. B. V., Amsterdam, 1995.

\bibitem[BM71]{Bruggesser:Mani}
H.~Bruggesser and P.~Mani.
\newblock Shellable decompositions of cells and spheres.
\newblock {\em Math. Scand.}, 29:197--205 (1972), 1971.

\bibitem[BW96]{Bjoerner:Wachs1}
A.~Bj{\"o}rner and M.~L. Wachs.
\newblock Shellable nonpure complexes and posets. {I}.
\newblock {\em Trans. Amer. Math. Soc.}, 348(4):1299--1327, 1996.

\bibitem[BW97]{Bjoerner:Wachs2}
A.~Bj{\"o}rner and M.~L. Wachs.
\newblock Shellable nonpure complexes and posets. {II}.
\newblock {\em Trans. Amer. Math. Soc.}, 349(10):3945--3975, 1997.

\bibitem[Can79]{Cannon}
J.~W. Cannon.
\newblock Shrinking cell-like decompositions of manifolds. {C}odimension three.
\newblock {\em Annals of Mathematics}, 110(1):83--112, 1979.

\bibitem[Cas95]{Cas}
M.~R. Casali.
\newblock A note about bistellar operations on {PL}-manifolds with boundary.
\newblock {\em Geom. Dedicata}, 56(3):257--262, 1995.

\bibitem[DK74]{Danaraj:Klee74}
G.~Danaraj and V.~Klee.
\newblock Shellings of spheres and polytopes.
\newblock {\em Duke Math. J.}, 41:443--451, 1974.

\bibitem[DK78]{Danaraj:Klee}
G.~Danaraj and V.~Klee.
\newblock A representation of {$2$}-dimensional pseudomanifolds and its use in
  the design of a linear-time shelling algorithm.
\newblock {\em Ann. Discrete Math.}, 2:53--63, 1978.
\newblock Algorithmic aspects of combinatorics (Conf., Vancouver Island, B.C.,
  1976).

\bibitem[dL13]{Longueville}
M.~de~Longueville.
\newblock {\em A course in topological combinatorics}.
\newblock Universitext. Springer, New York, 2013.

\bibitem[Edw75]{Edwards}
R.~D. Edwards.
\newblock The double suspension of a certain homology 3-sphere is {$S^5$}.
\newblock {\em Notices AMS}, 22(A-334), 1975.

\bibitem[Fis77a]{Fisk1977}
S.~Fisk.
\newblock Geometric coloring theory.
\newblock {\em Advances in Math.}, 24(3):298--340, 1977.

\bibitem[Fis77b]{Fisk}
S.~Fisk.
\newblock Variations on coloring, surfaces and higher-dimensional manifolds.
\newblock {\em Advances in Math.}, 25(3):226--266, 1977.

\bibitem[Hal04]{Hall}
H.~T. Hall.
\newblock {\em Counterexamples in discrete geometry}.
\newblock ProQuest LLC, Ann Arbor, MI, 2004.
\newblock Thesis (Ph.D.)--University of California, Berkeley.

\bibitem[HZ00]{Hachimori:Ziegler}
M.~Hachimori and G.~M. Ziegler.
\newblock Decompositons of simplicial balls and spheres with knots consisting
  of few edges.
\newblock {\em Math. Z.}, 235(1):159--171, 2000.

\bibitem[IJ03]{Izmestiev:Joswig}
I.~Izemestiev and M.~Joswig.
\newblock Branched coverings, triangulations, and 3-manifolds.
\newblock {\em Adv. Geom.}, 3(2):191--225, 2003.

\bibitem[IKN17]{IKN}
I.~Izmestiev, S.~Klee, and I.~Novik.
\newblock Simplicial moves on balanced complexes.
\newblock {\em Adv. Math.}, 320:82--114, 2017.

\bibitem[JKM17]{Juhnke:Murai}
M.~Juhnke-Kubitzke and S.~Murai.
\newblock {Balanced generalized lower bound inequality for simplicial
  polytopes}.
\newblock {\em Sel. Math. New Ser.}, pages 1--13, 2017.

\bibitem[KN16]{KN}
S.~Klee and I.~Novik.
\newblock Lower bound theorems and a generalized lower bound conjecture for
  balanced simplicial complexes.
\newblock {\em Mathematika}, 62(2):441--477, 2016.

\bibitem[Lic91]{Lickorish91}
W.~B.~R. Lickorish.
\newblock Unshellable triangulations of spheres.
\newblock {\em European J. Combin.}, 12(6):527--530, 1991.

\bibitem[Lic99]{Lickorish}
W.~B.~R. Lickorish.
\newblock Simplicial moves on complexes and manifolds.
\newblock 2:299--320, 11 1999.

\bibitem[LN16]{Lutz:Nevo}
F.~H. Lutz and E.~Nevo.
\newblock Stellar theory for flag complexes.
\newblock {\em Math. Scand.}, 118(1):70--82, 2016.

\bibitem[LR06]{Ludwig:Reitzner}
M.~Ludwig and M.~Reitzner.
\newblock Elementary moves on triangulations.
\newblock {\em Discrete Comput. Geom.}, 35:527--536, 2006.

\bibitem[Man16]{Manolescu}
C.~Manolescu.
\newblock Pin(2)-equivariant {S}eiberg-{W}itten {F}loer homology and the
  triangulation conjecture.
\newblock {\em J. Amer. Math. Soc.}, 29:147--176, 2016.

\bibitem[McM70]{McMullen}
P.~McMullen.
\newblock The maximum numbers of faces of a convex polytope.
\newblock {\em Mathematika}, 17:179--184, 1970.

\bibitem[MJKS17]{Juhnke:Murai:Novik:Sawaske}
I.~Novik M.~Juhnke-Kubitzke, S.~Murai and C.~Sawaske.
\newblock A generalized lower bound theorem for balanced manifolds.
\newblock {\em Math. Z.}, pages 1--22, 2017.

\bibitem[MS18]{Murai:Suzuki}
S.~Murai and Y.~Suzuki.
\newblock Balanced subdivisions and flips on surfaces.
\newblock {\em Proc. Amer. Math. Soc.}, 146:939--951, 2018.

\bibitem[New31]{Newman}
M.~H.~A. Newman.
\newblock A theorem in combinatory topology.
\newblock {\em J. London Math. Soc.}, S1-6(3):186--192, 1931.

\bibitem[Pac78]{Pac3}
U.~Pachner.
\newblock Bistellare \"aquivalenz kombinatorischer {M}annigfaltigkeiten.
\newblock {\em Arch. Math. (Basel)}, 30(1):89--98, 1978.

\bibitem[Pac91]{Pac1}
U.~Pachner.
\newblock P.{L}. homeomorphic manifolds are equivalent by elementary shellings.
\newblock {\em European J. Combin.}, 12(2):129--145, 1991.

\bibitem[RS82]{Rourke:Sanderson}
C.~P. Rourke and B.~J. Sanderson.
\newblock {\em Introduction to piecewise-linear topology}.
\newblock Springer Study Edition. Springer-Verlag, Berlin-New York, 1982.
\newblock Reprint.

\bibitem[Rud58]{Rudin}
M.~E. Rudin.
\newblock An unshellable triangulation of a tetrahedron.
\newblock {\em Bull. Amer. Math. Soc.}, 64:90--91, 1958.

\bibitem[Sta79]{St79}
R.~P. Stanley.
\newblock Balanced {C}ohen-{M}acaulay complexes.
\newblock {\em Trans. Amer. Math. Soc.}, 249(1):139--157, 1979.

\bibitem[Sta96]{Stanley-greenBook}
R.~P. Stanley.
\newblock {\em Combinatorics and commutative algebra}, volume~41 of {\em
  Progress in Mathematics}.
\newblock Birkh\"auser Boston, Inc., Boston, MA, second edition, 1996.

\bibitem[Wac07]{Wachs}
M.~L. Wachs.
\newblock Poset topology: tools and applications.
\newblock In {\em Geometric combinatorics}, volume~13 of {\em IAS/Park City
  Math. Ser.}, pages 497--615. Amer. Math. Soc., Providence, RI, 2007.

\bibitem[Zie95]{ZieglerBook}
G.~M. Ziegler.
\newblock {\em Lectures on Polytopes}.
\newblock Graduate texts in mathematics. Springer-Verlag, 1995.

\bibitem[Zie98]{Ziegler}
G.~M. Ziegler.
\newblock Shelling polyhedral {$3$}-balls and {$4$}-polytopes.
\newblock {\em Discrete Comput. Geom.}, 19(2):159--174, 1998.

\end{thebibliography}

\end{document}